\documentclass[11pt]{article}
\usepackage{amsmath,amssymb,amsthm,mathtools}
\usepackage[a4paper,margin=1in]{geometry}
\usepackage[T1]{fontenc}
\usepackage{lmodern}
\usepackage{microtype}
\usepackage[hidelinks]{hyperref}
\newtheorem{theorem}{Theorem}[section]
\newtheorem{lemma}[theorem]{Lemma}
\newtheorem*{lemma*}{Lemma}
\newtheorem{proposition}[theorem]{Proposition}
\newtheorem{corollary}[theorem]{Corollary}

\DeclareMathOperator{\ord}{ord}
\DeclareMathOperator{\rad}{rad}
\newcommand{\PP}{\mathcal P_{\mathrm{pow}}}

\hypersetup{pdftitle={Primary components and large prime factors of cyclotomic values},
pdfauthor={Wei Xie}}
\title{Primary components and large prime factors of cyclotomic values}
\author{Wei Xie}
\date{}
\begin{document}
\maketitle
\begin{abstract}
Let $\Phi_r(a,b)$ denote the homogeneous cyclotomic value.
For sufficiently large $A$, uniformly for $1/2\leq B\leq A$,
$r\geq3$ and $Z\geq2$, the number of coprime pairs in
$(A,2A]\times(B,2B]$ whose greatest prime factor $P(\Phi_r(a,b))$ is at most $Z$
is $O(AB/\varphi(r)+\sqrt{r/\varphi(r)}\,\pi(Z;r,1))$,
where $\pi(Z;r,1)$ counts primes at most $Z$ congruent to $1$ modulo $r$.
If $n=n(A)\to\infty$ through positive integers and
$1\leq H(A)=o(\log A)$, then $P(a^{2^n}+1)>2^nAH(A)$ for almost all $a\in(A,2A]$,
and $P(a^{2^n}+b^{2^n})>2^nA^2H(A)$ for almost all coprime pairs
in $(A,2A]^2$, with no further restriction on $n$.
For every fixed $\alpha>0$ and every $r\geq A^\alpha$,
pairwise distinct primitive prime divisors of $\Phi_r(a)$ can be
selected for all but $O_\alpha(A^{1-c_\alpha})$ bases in $(A,2A]$,
where $c_\alpha>0$.
For coprime positive integers $a,b$ with $\max(a,b)\geq2$
and $r\geq3$, put $C=\Phi_r(a,b)$ and
$S=\varphi(r)^2\log\max(a,b)$. Uniformly as $S\to\infty$,
$\max_{p\mid C}p^{v_p(C)}\geq(5/16-o(1))S$.
For every positive function $h(S)\to\infty$,
$\sum_{p\mid C,\;p^{v_p(C)}>S/h(S)}v_p(C)\log p
=(1-o(1))\log C$ uniformly in the same limit.
\end{abstract}

\section{Introduction}
Prime divisors of cyclotomic values are constrained by
multiplicative order. For positive coprime integers $a,b$ with
$\max(a,b)\geq2$ and $r\geq3$, put
\[
C_r(a,b)=\Phi_r(a,b)=b^{\varphi(r)}\Phi_r(a/b),
\]
where $\Phi_r$ is the $r$th cyclotomic polynomial and
$\Phi_r(a,b)$ denotes its homogeneous value.
Apart from a possible occurrence of the greatest prime divisor
of $r$, with multiplicity one, a prime divisor $p$ satisfies
$\ord_p(ab^{-1})=r$ and hence $p\equiv1\pmod r$.
The values grow exponentially with $\varphi(r)$, while their
prime factors lie in a progressively sparser residue class.

For an integer $N>1$, write $P(N)$ for its greatest prime factor
and $\PP(N)=\max_{p\mid N}p^{v_p(N)}$ for its largest primary
component. The factorization identity
\[
\log N=\sum_{p\mid N}v_p(N)\log p
\]
involves both the prime divisors and their multiplicities.
A large primary component may be a high power of a small
prime, so bounds for $\PP(N)$ and $P(N)$ are distinct questions.
Examples include the Fermat values
\[
\begin{aligned}
F_{a,n}&=a^{2^n}+1=C_{2^{n+1}}(a,1),\\
F_{a,b,n}&=a^{2^n}+b^{2^n}=C_{2^{n+1}}(a,b),\qquad n\geq1.
\end{aligned}
\]

Stewart~\cite{Stewart1975,Stewart1977}
obtained bounds for divisors of $a^r-b^r$ and recurrence sequences.
In particular, for fixed $a$, if $\Phi_r(a)$ is squarefree, this gives
$P(\Phi_r(a))\gg_a\varphi(r)^2$.
Quadratic lower bounds for primitive prime-power divisors of
Lucas numbers also originate in this work; Kiss~\cite[p.~775]{Kiss}
records a bound of order $\varphi(r)^2$.
The Lucas-sequence setting includes fixed coprime binary bases.
Murty and S\'eguin~\cite[Theorem~1.1]{MS} proved that, for fixed $a\geq2$,
if $v_p(a^{\ord_p(a)}-1)\leq\kappa$ for every prime $p\nmid a$,
then $P(\Phi_r(a))\gg_{a,\kappa}\varphi(r)^2$.
Xie~\cite{Xie} obtained
explicit squarefree estimates of order $4^n/n$.

Stewart's unconditional estimate~\cite[Theorem~1.1]{StewartLucasLehmer}
gives, for each fixed $a\geq2$ and all sufficiently large $n$,
\[
P(F_{a,n})>2^{n+1}
\exp\left(\frac{(n+1)\log2}{104\log((n+1)\log2)}\right).
\]
Bilu, Gun and Hong~\cite[Theorem~1.2 and (9.1)]{BGH}
proved $P(C_r(a,b))\geq r\exp(0.0005\log r/\log\log r)$
uniformly in the bases for $r\geq\exp(10^6)$.
For the classical Fermat numbers, further greatest-prime-factor
bounds were obtained by Le~\cite{Le}, Grytczuk, W\'ojtowicz and
Luca~\cite{GWL}, and Chen~\cite{Chen}; generalized Fermat numbers
were treated by Grytczuk and M\k{a}dryk~\cite{GM}.

When the bases vary, one can instead count values whose prime
divisors all lie below a given bound.
For a fixed irreducible binary form $F$ of degree $g\geq2$,
Balog, Blomer, Dartyge and Tenenbaum~\cite[Theorem~1]{BBDT}
proved that a positive proportion of the pairs $1\leq a,b\leq x$
give $y$-friable values when $y\geq x^{g-2+\varepsilon}$;
the constants depend on $F$ and on $\varepsilon>0$.
For the fixed cyclotomic index $r=12$,
Dartyge~\cite[Theorem~1.1]{Dartyge} proved a lower bound
$A^{1+c}$ for the greatest prime factor of
$\prod_{A<a\leq2A}\Phi_{12}(a)$, with explicit $c>0$.
In the family $C_r(a,b)$, varying $r$ changes both the degree
$\varphi(r)$ and the congruence condition $p\equiv1\pmod r$
for its primitive prime divisors.

We study lower bounds for $\PP(C_r(a,b))$ and $P(C_r(a,b))$
as the bases and the cyclotomic index vary.
For bases in an interval or a rectangle, we consider how often
$P(C_r(a,b))\leq Z$.
For the one-variable family $\Phi_r(a)$, we also ask whether
distinct bases can be assigned distinct primitive prime divisors.

\section{Preliminaries}
All logarithms are natural. Implied constants are absolute unless
dependence is indicated by a subscript. For a positive integer $N$,
write $\rad(N)=\prod_{p\mid N}p$. We write
$\pi(x;q,b)$ for the number of primes $p\leq x$ with
$p\equiv b\pmod q$, and
$\vartheta(x;q,b)=\sum_{p\leq x,\ p\equiv b\pmod q}\log p$.

\begin{lemma}[{Brun--Titchmarsh~\cite{MV,MontgomeryVaughan}; Iwaniec~\cite[Theorem~6]{Iw82}}]
\label{lem:classical-BT}\label{lem:iwaniec-BT}
Let $q\geq1$ be an integer and $(b,q)=1$.
The Brun--Titchmarsh inequality gives, for $x>q$,
$\pi(x;q,b)\leq2x/(\varphi(q)\log(x/q))$.
Iwaniec's bound gives, for each $\rho>0$, uniformly for $x^{2/5}<q\leq x^{2/3}$
and sufficiently large $x$,
\[
\pi(x;q,b)\leq
\frac{(8+\rho)x}{\varphi(q)\log(x^6q^{-7})}.
\]
\end{lemma}

\begin{lemma}\label{prop:BT}\label{lem:modulus-lifting}
For every $\varepsilon>0$, some $\delta>0$ and $x_0>0$ satisfy
\[
\vartheta(x;q,b)\leq(16/5+\varepsilon)x/\varphi(q)
\]
whenever $x\geq x_0$, $(b,q)=1$, and
$x^{1/2-\delta}\leq q\leq x^{1/2+\delta}$.
For $b=1$, the estimate also holds throughout
$2\leq q\leq x^{1/2+\delta}$.
\end{lemma}
\begin{proof}
Choose $\rho>0$ and $\delta<1/10$ so that
$(8+\rho)/(6-7(1/2+\delta))<16/5+\varepsilon$.
The second inequality of Lemma~\ref{lem:classical-BT}, followed
by $\vartheta\leq\pi\log x$, proves the first assertion.
For the second, choose the least $j\geq0$ with $Q=2^jq\geq\sqrt x$.
For large $x$, $\sqrt x\leq Q\leq\max(q,2\sqrt x)\leq x^{1/2+\delta}$.
Reduction of units modulo $Q$ onto units modulo $q$ is surjective,
with $\varphi(Q)/\varphi(q)$ classes above $1$.
Every prime $1\pmod q$ belongs to these classes: the only possible
new prime divisor of $Q$ is $2$, which is not $1\pmod q$.
Sum the first estimate over the classes.
\end{proof}

\begin{lemma}\label{primary:cyclotomic-facts}\label{lem:order-general}
\label{lem:pk-solutions}\label{lem:hom-facts}
Let $a,b$ be coprime positive integers with $\max(a,b)\geq2$,
and let $r\geq3$. A prime divisor of $C_r(a,b)$ either satisfies
$\ord_p(ab^{-1})=r$ and $p\equiv1\pmod r$, or is $P(r)$ and
occurs once. No prime divisor divides $ab$.
Uniformly for $a\geq2$ and $r\geq2$,
\[
\log\Phi_r(a)=\varphi(r)\log a+O(1/a).
\]
If $A<a,b\leq2A$, then
\begin{equation}\label{eq:hom-box}
A^{\varphi(r)}\Phi_r(1)\leq C_r(a,b)\leq(4A)^{\varphi(r)}.
\end{equation}
If $p\equiv1\pmod r$, then $\Phi_r(x)\equiv0\pmod{p^j}$
has exactly $\varphi(r)$ roots for every $j\geq1$; these are the
primitive $r$th roots in the $r$-element group of $r$th roots
modulo $p^j$. For a power of two $m\geq1$ and odd $p$, the equation
$x^m\equiv-1\pmod{p^j}$ has exactly $m$ roots if
$p\equiv1\pmod{2m}$ and none otherwise. When nonempty, these
roots form a coset of a subgroup of order $m$.
For every integer $h\geq1$,
\[
\Phi_{rh}(X)\mid\Phi_r(X^h)\quad\text{in }\mathbb Z[X],
\qquad C_{rh}(u,v)\mid C_r(u^h,v^h)
\]
for positive integers $u,v$. If $\rad(h)\mid r$, then
\[
\Phi_r(X^h)=\Phi_{rh}(X),\qquad \varphi(rh)=h\varphi(r).
\]
\end{lemma}
\begin{proof}
A prime dividing $a$ or $b$ does not divide $C_r(a,b)$,
since the constant and leading coefficients of $\Phi_r$ are $1$.
For odd $p\nmid ab$, put $d=\ord_p(ab^{-1})$ and
$w=v_p(a^d-b^d)$. If $d\nmid r$, the valuation is zero.
Factoring a difference of powers shows that an exponent coprime
to $p$ leaves the valuation unchanged. By the binomial theorem,
multiplying the exponent by $p$ adds one. Hence
$v_p(a^{df}-b^{df})=w+v_p(f)$ for $f\geq1$.
From $a^r-b^r=\prod_{e\mid r}\Phi_e(a,b)$,
M\"obius inversion gives, for $r=dt$,
\[
v_p(C_r(a,b))=\sum_{f\mid t}\mu(t/f)(w+v_p(f)).
\]
Separating the $p$-part of $t$, this is $w$ for $t=1$, is $1$
when $t$ is a positive power of $p$, and is zero otherwise. In the second case
$r=dp^j$ with $d\mid p-1$, so $p=P(r)$.
When $a,b$ are odd, the usual formula
$v_2(a^e-b^e)=v_2(a-b)+v_2(a+b)+v_2(e)-1$ for even $e$,
and $v_2(a^e-b^e)=v_2(a-b)$ for odd $e$, gives
$v_2(C_r(a,b))=1$ for $r=2^j$, $j\geq2$, and zero for
the other $r\geq3$.
The product formula for $\Phi_r$ bounds the unary logarithmic
error by
\[
\sum_{d\geq1}-\log(1-a^{-d})
\leq\sum_{d\geq1}\frac{a^{-d}}{1-a^{-d}}
\leq\frac{a}{(a-1)^2}=O(1/a).
\]
For a primitive $r$th root of unity $\zeta$,
\[
|a-b\zeta|^2=(a-b)^2+ab|1-\zeta|^2
\geq A^2|1-\zeta|^2.
\]
Multiplication over the roots gives \eqref{eq:hom-box};
$\Phi_r(1)\geq1$, and $|a-b\zeta|\leq4A$ gives the upper bound.
If $p\nmid r$, each $r$th root lifts uniquely from $p^j$ to $p^{j+1}$:
in
\[
(x+tp^j)^r\equiv x^r+rx^{r-1}tp^j\pmod{p^{j+1}},
\]
the coefficient $rx^{r-1}$ is a unit modulo $p$.
When $p\equiv1\pmod r$, the $r$th roots modulo $p$ form a cyclic
group of order $r$, with $\varphi(r)$ primitive roots.
Reduction identifies the root groups at all levels; their primitive
roots are exactly the roots of $\Phi_r$ by the valuation formula above.
For a power of two $m$, the roots of $x^m=-1$ are precisely the
primitive $2m$th roots. For odd $p$, they exist exactly when
$2m\mid p-1$, number $m$ at every level, and form a coset of
the subgroup of $m$th roots.
Finally, the $h$th power of a primitive $(rh)$th root has
order $r$, so $\Phi_{rh}(X)$ divides $\Phi_r(X^h)$.
Division by a monic integer polynomial gives an integer
quotient, and homogenization gives the assertion for $C_r$.
If $\rad(h)\mid r$, the totient product formula gives
$\varphi(rh)=h\varphi(r)$; the preceding divisibility is then
an equality of monic polynomials of the same degree.
\end{proof}

\begin{lemma}\label{bi-fractions}\label{lem:fractions}
Let $M\geq2$, $U,V\geq1$, and
$\lambda\in(\mathbb Z/M\mathbb Z)^\times$.
The number of reduced positive fractions $u/v$ with
$u\leq U$, $v\leq V$, $\gcd(uv,M)=1$, and
$u\equiv\lambda v\pmod M$ is at most
\[
2(1+UV/M).
\]
If $M>UV$, reduction modulo $M$ is injective on these
fractions. In particular, for $H\geq1$ and $u,v\leq H$ the bound is
$2(1+H^2/M)$, and $M>H^2$ gives injectivity.
\end{lemma}
\begin{proof}
Two distinct fractions in one class satisfy
$|uv'-u'v|\geq M$. Those at most $U/V$ lie in
$(0,U/V]$ and are separated by at least $M/V^2$, so there
are at most $1+UV/M$ of them. For the others, take inverses:
they lie in $(0,V/U)$, in the class $\lambda^{-1}$, and
are separated by at least $M/U^2$.
If $M>UV$, a nonzero cross-product difference has absolute
value less than $UV<M$ and cannot be divisible by $M$.
\end{proof}

\begin{lemma}\label{lem:dependent}\label{bi-dependent}
Let $A>1$. Suppose that $a/b$ and $c/d$ are reduced
positive fractions with
$A<\max(a,b),\max(c,d)\leq2A$.
If they are multiplicatively dependent, then
$(c,d)=(a,b)$ or $(c,d)=(b,a)$.
In particular, distinct integers in $(A,2A]$ are
multiplicatively independent.
\end{lemma}
\begin{proof}
Dependent fractions have the form $(u/v)^s$ and $(u/v)^t$,
where $u>v\geq1$, $(u,v)=1$, and $s,t$ are nonzero integers.
This follows by dividing their common prime-exponent direction
by the greatest common divisor of its coordinates.
Their maximum numerator or denominator is $u^{|s|}$ and
$u^{|t|}$. If $|s|\ne|t|$, their ratio is at least
$u\geq2$, contradicting the window. Thus $t=\pm s$.
Reduction of the fractions gives the identical or reversed
pairs. For distinct integers in $(A,2A]$, neither equality
nor inversion is possible.
\end{proof}

\begin{lemma}\label{lem:progression}\label{bi-progression}
Uniformly for integers $r\geq2$ and $X\geq2$,
\begin{equation}\label{eq:primecount}
\pi(X;r,1)\ll\min\left\{\frac Xr,
\frac{X}{\varphi(r)\log(2+X/r)}\right\},\qquad
\vartheta(X;r,1)\leq\pi(X;r,1)\log X.
\end{equation}
In particular, for powers of two $m\geq2$,
\[
\pi(X;2m,1)\ll\frac{X}{m\log(2+X/m)},\qquad
\vartheta(X;2m,1)\ll\frac{X\log X}{m\log(2+X/m)}.
\]
\end{lemma}
\begin{proof}
Comparison with the integers $1+\ell r$, $\ell\geq1$,
gives the first bound in the minimum. For $X\leq2r$ it
also gives the second, since $\varphi(r)\leq r$.
For $X>2r$, Brun--Titchmarsh gives the second bound.
Multiplication by $\log X$ gives the estimate for
$\vartheta$, and $\varphi(2m)=m$ gives the special case.
\end{proof}

For a fixed integer $s\geq2$, let $d_s(u)$ be the number of
ordered representations of a positive integer $u$ as a product
of $s$ positive integers. Unique factorization gives
\[
d_s(u)=\prod_{\ell^\nu\parallel u}
\binom{\nu+s-1}{s-1},
\qquad d_s(u)\ll_{s,\eta}u^\eta\quad(\eta>0).
\]
Indeed, $\binom{\nu+s-1}{s-1}\leq s^\nu$ for $\nu\geq1$,
so the factor at every sufficiently large prime $\ell$ is at
most $\ell^{\eta\nu}$. For each of the finitely many remaining
primes, its ratio to $\ell^{\eta\nu}$ is bounded, since the
numerator grows polynomially in $\nu$.

\begin{lemma}\label{lem:amplification}
Let $M\geq2$, let $H$ be a subgroup of
$(\mathbb Z/M\mathbb Z)^\times$ of order $t$, and let $C$ be
a coset of $H$. For $A\geq2$, let $\mathcal R$ be any set of
integers in $(A,2A]$ whose residues lie in $C$.
For an integer $s\geq2$, if every integer
has at most $D_s$ ordered representations as a product of $s$ members of
$\mathcal R$, then
\begin{equation}\label{eq:prodcount}
|\mathcal R|^s\leq D_s t\bigl((2A)^s/M+1\bigr).
\end{equation}
Consequently, for every fixed integer $s\geq2$ and every
$\varepsilon>0$, uniformly for $A\geq2$,
\[
\#\{A<a\leq2A:a\bmod M\in C\}
\ll_{s,\varepsilon}A^\varepsilon
\left\{A(t/M)^{1/s}+t^{1/s}\right\}.
\]
\end{lemma}
\begin{proof}
The $|\mathcal R|^s$ ordered products are at most $(2A)^s$, and
their residues lie in $C^s$, which has $t$ elements.
Each residue class contains at most $(2A)^s/M+1$ possible
integer products, and each product has at most $D_s$
representations. This proves \eqref{eq:prodcount}.
For the last assertion, take all interval integers in the
coset and use $d_s(u)\ll_{s,\varepsilon}u^\varepsilon$,
so that $D_s\ll_{s,\varepsilon}A^{s\varepsilon}$.
Taking the $s$th root gives the stated bound.
\end{proof}

In particular, let $m\geq2$ be a power of two and let
$p\equiv1\pmod{2m}$ be an odd prime. The roots of
$x^m\equiv-1\pmod{p^j}$ form a coset of a subgroup of order
$m$, by Lemma~\ref{lem:pk-solutions}. Thus, for every
$j\geq1$,
\begin{equation}\label{eq:third-root}
\#\{A<a\leq2A:p^j\mid a^m+1\}
\ll_\varepsilon A^\varepsilon
\left(Am^{1/3}p^{-j/3}+m^{1/3}\right).
\end{equation}

\begin{lemma}[{Matolcsi--Ruzsa~\cite[Theorem~1.5]{MR}}]\label{lem:planar}\label{bi-planar}
Let $X\subset\mathbb R^2$ be a finite set not contained
in a line. For each integer $h\geq1$, let $hX$ denote
the set of sums of $h$ members of $X$, with repetitions
allowed. Then
\[
|hX|\geq(|X|-2)\frac{h(h+1)}2+h+1.
\]
\end{lemma}
\begin{proof}
For $h=1$ the bound is an equality. Let $h\geq2$, and
triangulate $\operatorname{conv}(X)$ with vertex set $X$.
Euler's formula, together with counting the edges of the
triangles, gives at least $|X|-2$ triangles and $2|X|-3$ edges.
After scaling the triangulation by $h$, each vertex belongs
to $hX$, each open edge contributes $h-1$ points of $hX$,
and each open triangle contains at least
$(h-1)(h-2)/2$ points of $hX$, obtained from positive
integer barycentric coordinates with sum $h$.
These relative interiors are disjoint. Hence
\[
\begin{aligned}
|hX|&\geq |X|+(2|X|-3)(h-1)
       +(|X|-2)\frac{(h-1)(h-2)}2\\
&=(|X|-2)\frac{h(h+1)}2+h+1.
\end{aligned}
\]
\end{proof}

\begin{lemma}\label{lem:growth}\label{hom-growth}
Let $A>1$ and let $h\geq1$ be an integer. Let $\mathcal R$
consist either of $N\geq1$ distinct integers in $(A,2A]$, or
of $N\geq1$ distinct reduced fractions $a/b$ with
$A<a\leq2A$ and $B<b\leq2B$, where $1/2\leq B\leq A$.
In the integer case, the products of at most $h$ members,
with repetitions allowed and including the empty product,
number at least
\begin{equation}\label{eq:growth}
(N-1)\frac{h(h+1)}2+h+1.
\end{equation}
The same bound holds in the rational case when $N\geq3$.
Put $\kappa=1$ in the integer case and $\kappa=2$ in the
rational case. If the reduction of $\mathcal R$ modulo
$p^j$ lies in a multiplicative group of order $r$, then
\begin{equation}\label{eq:root-growth}
N\leq\max\left\{\kappa,1+\frac{2r}{h(h+1)}\right\}
\end{equation}
provided that $p^j>(2A)^h$ in the integer case and
$p^j>(4AB)^h$ in the rational case. Here all numerators
and denominators are assumed coprime to $p$.
\end{lemma}
\begin{proof}
For one integer the empty product and its first $h$ powers give
$h+1$ products. Otherwise, in the integer case $N\geq2$, and
in the rational case we need only consider $N\geq3$.
Every member of $\mathcal R$ differs from $1$. Also
$\max(a,b)\in(A,2A]$ for the reduced fractions.
Lemma~\ref{lem:dependent} gives multiplicatively independent
$x,y\in\mathcal R$.

Choose a prime $\ell$ with $v_\ell(x)\ne0$, where
$v_\ell(a/b)=v_\ell(a)-v_\ell(b)$, and put
\[
X=\{(0,0)\}\cup
\{(\log z,v_\ell(z)):z\in\mathcal R\}.
\]
The first coordinate gives $|X|=N+1$. If the points corresponding
to $x,y$ were collinear with the origin, then
\[
v_\ell(x)\log y=v_\ell(y)\log x,
\qquad y^{v_\ell(x)}=x^{v_\ell(y)},
\]
contrary to multiplicative independence. Thus $X$ is not all
on one line. The map $z\mapsto(\log z,v_\ell(z))$ is injective
and sends multiplication to addition. Consequently, the products
of at most $h$ members of $\mathcal R$ correspond bijectively
to $hX$, with the origin padding shorter products.
Lemma~\ref{lem:planar} now gives \eqref{eq:growth}.

Integer products have height at most $(2A)^h$, so reduction
modulo $p^j$ is injective when $p^j>(2A)^h$. Rational
products have reduced numerator at most $(2A)^h$ and
denominator at most $(2B)^h$; these bounds also hold for
shorter products because $2B\geq1$. Lemma~\ref{lem:fractions}
with $U=(2A)^h$ and $V=(2B)^h$ gives injectivity when
$p^j>(4AB)^h$.
All product residues remain in the same group. Its order
therefore bounds \eqref{eq:growth} from above when
$N>\kappa$, giving \eqref{eq:root-growth}. For
$N\leq\kappa$ the conclusion holds directly.
\end{proof}

\begin{corollary}\label{lem:pairmass}
Let $A\geq2$, let $a,b$ be distinct integers in $(A,2A]$,
and let $r\geq3$. If $p\equiv1\pmod r$ is prime and
$p^j$ divides both $\Phi_r(a)$ and $\Phi_r(b)$ for an
integer $j\geq1$, then, with $s=\gcd(r,2)$,
\[
j\log p\leq2\sqrt{\frac{r}{3s}}\log(2A)+\log s.
\]
\end{corollary}
\begin{proof}
By Lemma~\ref{primary:cyclotomic-facts}, the residues of
$a,b$ generate the $r$-element group $G$ of $r$th roots
modulo $p^j$. Put $E=G\cap\{1,-1\}$, so $|E|=s$.
The kernel $\Lambda$ of the surjective homomorphism
\[
\mathbb Z^2\longrightarrow G/E,\qquad
(u,v)\longmapsto a^ub^vE,
\]
is a lattice of determinant $r/s$. The centrally symmetric
hexagon
\[
B_H=\{(u,v)\in\mathbb R^2:
\max(|u|,|v|,|u+v|)\leq H\}
\]
has area $3H^2$. If $H>2\sqrt{r/(3s)}$, the area of
$B_H/2$ exceeds $\det\Lambda$. Averaging
$\#((x+\Lambda)\cap(B_H/2))$ over a fundamental
parallelogram gives a translate containing two points.
Their difference is a nonzero lattice point in $B_H$.
Letting $H$ decrease to the endpoint, and using finiteness
of lattice points in a fixed bounded set, gives
$(u,v)\in\Lambda\setminus\{(0,0)\}$ with
\[
\max(u_++v_+,u_-+v_-)
=\max(|u|,|v|,|u+v|)\leq2\sqrt{r/(3s)},
\]
where $u_\pm=\max(\pm u,0)$, and similarly for $v$.
For some $\epsilon\in E$, the integer
\[
D=a^{u_+}b^{v_+}-\epsilon a^{u_-}b^{v_-}
\]
is divisible by $p^j$. It is nonzero: for $\epsilon=1$,
use Lemma~\ref{lem:dependent}; for $\epsilon=-1$, it is
a sum of positive integers. Finally,
$p^j\leq|D|\leq s(2A)^{2\sqrt{r/(3s)}}$, which proves
the bound.
\end{proof}

For real $X>r\geq3$, put
\begin{equation}\label{eq:first-moment-H}
\mathcal H(r,X)=\log(2r)\log\bigl(2+\log(X/r)\bigr)+\log(X/r).
\end{equation}

\begin{lemma}\label{lem:rpharm}
Uniformly for $X>r\geq3$,
\[
\sum_{\substack{p\leq X\\p\equiv1\pmod r}}\frac{\log p}{p}
\ll\frac{\mathcal H(r,X)}{\varphi(r)}.
\]
Consequently, uniformly for $X\geq3$ and integers $r\geq3$,
the same sum is $O(\log X\log(2+\log X)/\varphi(r))$.
\end{lemma}
\begin{proof}
There is at most one possible prime in
$(r,2r]$, contributing $O(\log(2r)/\varphi(r))$.
This proves the first bound for $r<X\leq2r$.
For $X>2r$, Brun--Titchmarsh and partial summation bound
the remaining sum by
\[
\frac{2\log X}{\varphi(r)\log(X/r)}+
\frac2{\varphi(r)}\int_{2r}^X\frac{\log t}{t\log(t/r)}\,dt.
\]
Writing $u=\log(X/r)>\log2$, the integral equals
\[
u-\log2+\log r\log\frac{u}{\log2}.
\]
It is $O(\mathcal H(r,X))$, as is the endpoint, since
$\log X/\log(X/r)=1+\log r/u\leq1+\log r/\log2$.
The second assertion follows from
$\mathcal H(r,X)\ll\log X\log(2+\log X)$ when $X>r$;
the sum is empty otherwise.
\end{proof}

% Merger reconstruction, v24. The source manuscripts remain unchanged.
% Shared inputs: primary:cyclotomic-facts (including primitive binary pairs),
% prop:BT, lem:modulus-lifting, and the totient lower bound cited as RS.

\section{Primary components and radical estimates}\label{sec:primary}

For positive coprime integers \(a,b\) with \(B=\max(a,b)\geq2\),
and for \(r\geq3\), put
\[
C_r(a,b)=b^{\varphi(r)}\Phi_r(a/b),\qquad
q=\varphi(r),\qquad S=q^2\log B.
\]
Recall that \(\PP(N)=\max_{p\mid N}p^{v_p(N)}\).
The estimates below include \(\Phi_r(a)=C_r(a,1)\).

The index $r$ and the degree $q$ satisfy the classical
bound~\cite[Theorem~15]{RS}
\begin{equation}\label{eq:totient}
\frac rq\ll\log\log(3r)\qquad(r\geq3).
\end{equation}

\subsection{Primary components}

\begin{lemma}\label{lem:hom-size}
Uniformly for the preceding parameters,
\[
\log C_r(a,b)=q\log B+O(\mathrm d(r)\log r),
\]
where \(\mathrm d(r)\) is the number of positive divisors of \(r\).
Consequently,
\[
\log C_r(a,b)=(1+o(1))q\log B
\qquad(S\to\infty)
\]
uniformly in \(a,b,r\).
\end{lemma}
\begin{proof}
Set \(\lambda=\min(a,b)/B\). Since \((a,b)=1\) and
\(B\geq2\), we have \(0<\lambda<1\). Reciprocity of \(\Phi_r\) gives
\(C_r(a,b)=B^q\Phi_r(\lambda)\). Since
\(\sum_{d\mid r}\mu(r/d)=0\), the product formula yields
\[
\log\Phi_r(\lambda)
=\sum_{d\mid r}\mu(r/d)
\log\frac{1-\lambda^d}{1-\lambda}.
\]
For every \(d\geq1\), the quotient is
\(1+\lambda+\cdots+\lambda^{d-1}\), between \(1\) and \(d\).
Thus the absolute value of the sum is at most
\(\mathrm d(r)\log r\).
The elementary bound \(\mathrm d(r)\leq2\sqrt r\)
and \eqref{eq:totient} show that
\(\mathrm d(r)\log r=o(q)\) as \(r\to\infty\).
Since \(\log B\geq\log2\), this proves the relative estimate
for large indices. The finitely many remaining indices are
handled by \(\log B\to\infty\).
\end{proof}

\begin{theorem}\label{thm:mass}
\label{primary:cyclotomic-extension}\label{primary:mass}
Fix \(c>0\). Uniformly for positive coprime integers \(a,b\)
with \(B\geq2\) and integers \(r\geq3\), as \(S\to\infty\),
\[
\sum_{\substack{p\mid C_r(a,b)\\p^{v_p(C_r(a,b))}\leq cS}}
v_p(C_r(a,b))\log p
\leq\left(\frac{16c}{5}+o_c(1)\right)q\log B.
\]
Consequently,
\[
\PP(C_r(a,b))\geq\left(\frac5{16}-o(1)\right)S.
\]
For every positive function \(h(S)\to\infty\),
\[
\log\left(
\prod_{\substack{p\mid C_r(a,b)\\p^{v_p(C_r(a,b))}>S/h(S)}}
p^{v_p(C_r(a,b))}\right)
=(1-o(1))\log C_r(a,b).
\]
For each fixed integer \(k\geq1\), if all prime multiplicities
in \(C_r(a,b)\) are at most \(k\), then
\[
P(C_r(a,b))\geq\left(\frac5{16k}-o_k(1)\right)S.
\]
\end{theorem}

\begin{proof}
Write $L=\log B$, $C=C_r(a,b)$, and $e_p=v_p(C)$.
For every $Y\geq2$, Lemma~\ref{primary:cyclotomic-facts} gives
\begin{equation}\label{eq:lowmass}
\begin{split}
&\sum_{\substack{p\mid C\\p\mid r}}\log p+
\sum_{\substack{p\mid C\\p\equiv1\pmod r}}
\sum_{\substack{1\leq j\leq e_p\\p^j\leq Y}}\log p\\
&\qquad\leq\log r+\vartheta(Y;r,1)
+\frac{\sqrt Y}{r}\log Y.
\end{split}
\end{equation}
The first sum includes the possible exceptional prime to its
full multiplicity, which is one. For the ordinary primes,
the first layers contribute at most $\vartheta(Y;r,1)$.
An additional layer requires $p\leq\sqrt Y$; there are
at most $\sqrt Y/r$ such primes, and each contributes at
most $\log Y$ over all its layers. Every component
$p^{e_p}\leq Y$ is covered by the left-hand side.

Now put $Y=cS$. By \eqref{eq:totient},
\[
\frac r{\sqrt Y}=\frac{r}{\sqrt c\,q\sqrt L}
\ll_c\log\log(3r).
\]
For each fixed $\delta>0$, this gives $r\leq Y^{1/2+\delta}$
uniformly once $Y$ is sufficiently large. For large $r$,
use $Y\geq c(\log2)q^2$; bounded indices are handled by
$Y\to\infty$. Lemma~\ref{lem:modulus-lifting} yields
\[
\vartheta(Y;r,1)
\leq\left(\frac{16}{5}+o_c(1)\right)\frac Yq
=\left(\frac{16c}{5}+o_c(1)\right)qL.
\]
The other terms in \eqref{eq:lowmass} are $o_c(qL)$.
For the last term, $q\leq r$ gives
\[
\frac{\sqrt Y\log Y}{rqL}
\ll_c\frac{\log S}{\sqrt S}\longrightarrow0.
\]
Also $\log r=o(qL)$ uniformly: use $\log r=o(q)$ for
large $r$ and $L\to\infty$ for bounded $r$.
This proves the asserted small-component mass bound.

Lemma~\ref{lem:hom-size} gives
\(\log C_r(a,b)=(1+o(1))qL\) uniformly as \(S\to\infty\).
For each fixed \(c<5/16\), the mass bound therefore excludes
the possibility that all components are at most \(cS\).
Letting \(c\) approach \(5/16\) proves the bound for \(\PP\).
Fix \(c>0\). Eventually \(S/h(S)\leq cS\), so
\[
\limsup_{S\to\infty}
\frac{\displaystyle\sum_{\substack{p\mid C_r(a,b)\\
p^{v_p(C_r(a,b))}\leq S/h(S)}}v_p(C_r(a,b))\log p}
{\log C_r(a,b)}\leq\frac{16c}{5}.
\]
Letting \(c\to0\) proves the concentration assertion.
Finally, if every multiplicity is at most \(k\) and
\(P(C_r(a,b))\leq cS\), then
\[
\log C_r(a,b)\leq\log r+k\vartheta(cS;r,1)
\leq\left(\frac{16ck}{5}+o_{c,k}(1)\right)qL.
\]
This contradicts Lemma~\ref{lem:hom-size} for every fixed
\(c<5/(16k)\) once \(S\) is sufficiently large.
\end{proof}

If $a=u^h$, $b=v^h$
for coprime positive integers $u,v$ and an integer $h\geq1$,
then Lemma~\ref{primary:cyclotomic-facts} and the preceding
theorem give
\[
\PP(C_r(a,b))\geq\left(\frac5{16}-o(1)\right)
\frac{\varphi(rh)^2}{h}\log B,
\]
uniformly as $\varphi(rh)^2\log\max(u,v)\to\infty$.
When $\rad(h)\mid r$, we have $C_r(a,b)=C_{rh}(u,v)$ and
\[
\varphi(rh)^2\log\max(u,v)=hS,\qquad
\varphi(rh)\log\max(u,v)=q\log B.
\]
Theorem~\ref{thm:mass} therefore applies with $hS$ in place
of $S$, also for varying $h$ with $hS\to\infty$.

Taking $b=1$ gives all the preceding assertions for
$\Phi_r(a)$, with $S=\varphi(r)^2\log a$.

For integers \(n\geq1\), write
\(F_{a,n}=a^{2^n}+1=\Phi_{2^{n+1}}(a)\); then
\(\varphi(2^{n+1})=2^n\). In particular, uniformly as
\(4^n\log a\to\infty\),
\[
\PP(F_{a,n})\geq\left(\frac5{16}-o(1)\right)4^n\log a,
\]
and, for every fixed \(c>0\),
\[
\prod_{\substack{p\mid F_{a,n}\\
p^{v_p(F_{a,n})}>c4^n\log a}}p^{v_p(F_{a,n})}
\geq a^{2^n(1-16c/5-o_c(1))}.
\]
The latter follows by subtracting the mass bound from
\(\log F_{a,n}>2^n\log a\).
\label{primary:near-quadratic-mass}
For each fixed \(0<\varepsilon<1\),
\[
\sum_{\substack{p\mid F_{a,n}\\p^{v_p(F_{a,n})}>2^{(2-\varepsilon)n}}}
v_p(F_{a,n})\log p=(1-o(1))\log F_{a,n}
\]
uniformly for \(a\geq2\) as \(n\to\infty\), since
\(2^{-\varepsilon n}/\log a\to0\).
\label{primary:bounded-multiplicity}
For each fixed \(\ell\geq2\), an \(\ell\)-free Fermat value
satisfies
\[
P(F_{a,n})\geq
\left(\frac5{16(\ell-1)}-o_\ell(1)\right)4^n\log a
\]
uniformly as \(4^n\log a\to\infty\).

\begin{lemma}\label{lem:quadratic}\label{lem:quadratic-restriction}
Let $a,b$ be positive coprime integers with $\max(a,b)\geq2$,
and let $r\geq3$. Let $d$ be the squarefree part of $ab$, and put
\[
D=\begin{cases}d,&d\equiv1\pmod4,\\4d,&d\not\equiv1\pmod4,\end{cases}
\qquad Q=\operatorname{lcm}(2r,D).
\]
For $D=1$, let $\chi_D=1$; otherwise $\chi_D$ is the
quadratic character of conductor $D$, given at odd primes
not dividing $ab$ by $\chi_D(p)=(d/p)$.
Every nonexceptional prime divisor of $C_r(a,b)$ satisfies
\[
\chi_D(p)=(-1)^{(p-1)/r}.
\]
Among the reduced residue classes $t\pmod Q$ with
$t\equiv1\pmod r$, this restriction selects exactly
$\varphi(Q)/(2\varphi(r))$ classes, except in the following
cases, when it imposes no additional restriction:
\[
r\text{ odd and }D\mid r,
\qquad\text{or}\qquad
r\text{ even},\ D\mid2r,\ D\nmid r.
\]
In particular, the restriction always selects half the classes
when $8\mid r$. In that case, if $c$ is the odd squarefree
part of $ab$, it is equivalent to
$(p/c)=(-1)^{(p-1)/r}$, with $(t/1)=1$, and
$Q=\operatorname{lcm}(2r,c)$.
At $b=1$ and $r=2^{n+1}$, where $n\geq2$, the
number of admissible classes is $\varphi(c)$ modulo $2^{n+2}c$.
For $a=2$ and $n\geq2$, every prime divisor of $F_{2,n}$ is $1\pmod{2^{n+2}}$.
\end{lemma}

\begin{proof}
Quadratic reciprocity gives $\chi_D$ exact conductor $D$.
For $N\mid Q$, it is trivial on the kernel of unit reduction
modulo $N$ if $D\mid N$. Otherwise, choose a nonresidue at
an odd prime of $D$ absent from $N$, or choose $3\pmod4$ or
$5\pmod8$ when the required $2$-part is missing.
The other components can be chosen to be $1$ by the Chinese
remainder theorem. This gives a kernel element with character
$-1$, and proves the converse.

For a nonexceptional prime $p$, the order of $ab^{-1}$
modulo $p$ is $r$. Euler's criterion therefore gives
\[
\chi_D(p)=\left(\frac{ab}{p}\right)
=\left(\frac{ab^{-1}}p\right)=(-1)^{(p-1)/r}.
\]
On
\[
H=\{t\in(\mathbb Z/Q\mathbb Z)^\times:t\equiv1\pmod r\},
\]
put $\epsilon(t)=(-1)^{(t-1)/r}$.
This is a character. If $r$ is even, multiplicativity follows
by expanding $(1+ru)(1+rv)$; if $r$ is odd, every unit $t$
is odd, so $\epsilon=1$ on $H$. Reduction onto the units
modulo $r$ is surjective, and hence $|H|=\varphi(Q)/\varphi(r)$.

For odd $r$, the character $\chi_D\epsilon$ on $H$ is
trivial precisely when $D\mid r$, by the conductor property
just proved. For even $r$, suppose it is trivial. The kernel
of reduction modulo $2r$ lies in $H$ and has $\epsilon=1$,
so $D\mid2r$. On the other hand, $\epsilon$ is nontrivial:
the class $1+r$ modulo $2r$ is a unit and lifts to an element
of $H$. Thus $D\nmid r$. Conversely, if $D\mid2r$ and
$D\nmid r$, both $\chi_D$ and $\epsilon$ restrict to the nontrivial
character of the two-element kernel of
\[
(\mathbb Z/2r\mathbb Z)^\times
\longrightarrow(\mathbb Z/r\mathbb Z)^\times,
\]
and their product is trivial on $H$.
In every other case, $\chi_D\epsilon$ is nontrivial, so its
kernel has exactly half the elements of $H$.

Since $D$ has $2$-part at most $8$, the even-index exception
cannot occur when $8\mid r$. Then $t\equiv1\pmod8$, and
quadratic reciprocity reduces $\chi_D(t)$ to $(t/c)$.
For $r=2^{n+1}$ and $b=1$, we have $Q=2^{n+2}c$ and
$\varphi(Q)/(2\varphi(r))=\varphi(c)$.
For $a=2$, $c=1$ and the restriction is $p\equiv1\pmod{2^{n+2}}$;
there is no exceptional prime divisor.
\end{proof}

\begin{corollary}\label{cor:fixed}
\label{primary:cyclotomic-fixed-base}\label{primary:fixed-base}
Fix positive coprime integers $a,b$ with $B=\max(a,b)\geq2$,
and let $D$ be as in Lemma~\ref{lem:quadratic}.
As $r\to\infty$ through odd indices with $D\nmid r$, or
through even indices with $D\nmid2r$ or $D\mid r$,
\[
\PP(C_r(a,b))\geq
\left(\frac58-o_{a,b}(1)\right)\varphi(r)^2\log B.
\]
In particular, the estimate holds whenever $8\mid r$, and
$b=1$ gives the one-variable estimate.
If $a=u^h$, $b=v^h$, with fixed positive coprime integers
$u,v$ and a fixed integer $h\geq1$, then, as $r\to\infty$ with $8\mid rh$,
\[
\PP(C_r(a,b))\geq
\left(\frac58-o_{a,b,h}(1)\right)
\frac{\varphi(rh)^2}{h}\log B.
\]
When $\rad(h)\mid r$, this is
$(5h/8-o_{a,b,h}(1))\varphi(r)^2\log B$.
For a fixed integer base $a\geq2$, set
\[
\tau(a)=\max\{t\in\mathbb Z_{\geq0}:a=u^{2^t}
\text{ for an integer }u\geq2\},
\]
and put
\[
\gcd_{p\mid a}v_p(a)=2^{\tau(a)}k(a),\qquad
\kappa(a)=2^{\tau(a)}\frac{\varphi(k(a))^2}{k(a)},
\]
where $k(a)$ is odd. For every fixed base of either parity,
\[
\PP(F_{a,n})\geq
\left(\frac{5\kappa(a)}8-o_a(1)\right)4^n\log a
\qquad(n\to\infty).
\]
Here $\kappa(a)\geq2^{\tau(a)}$, with strict inequality
when $k(a)>1$.
\end{corollary}

\begin{proof}
Fix $c>0$ and put $C=C_r(a,b)$, $q=\varphi(r)$, and
$Y=cq^2\log B$.
Since $D$ is fixed, $Q/r\leq2D$, while
$\log\varphi(r)=\log r+o(\log r)$. Hence
$\log Q/\log Y=1/2+o_{a,b,c}(1)$.
Under the stated index conditions, Lemma~\ref{lem:quadratic}
gives $\varphi(Q)/(2\varphi(r))$ admissible classes.
Summing Lemma~\ref{prop:BT} over them yields
\[
\sum_{\substack{p\leq Y\\p\text{ admissible}}}\log p
\leq\left(\frac85+o_{a,b,c}(1)\right)
\frac{Y}{\varphi(r)}.
\]
All nonexceptional prime divisors lie in these classes.
The argument of \eqref{eq:lowmass}, with its first layers
restricted to the admissible classes, therefore gives
\[
\sum_{\substack{p\mid C\\p^{v_p(C)}\leq Y}}v_p(C)\log p
\leq\left(\frac{8c}{5}+o_{a,b,c}(1)\right)q\log B.
\]
Here the additional layers and the exceptional prime cost
$O((\sqrt Y/r)\log Y+\log r)=o_{a,b,c}(q\log B)$.
Lemma~\ref{lem:hom-size} therefore excludes the possibility
that all primary components are at most $c\varphi(r)^2\log B$
for each $c<5/8$, once $r$ is sufficiently large.
This proves the first assertion, including the stated
specializations.

For $a=u^h$ and $b=v^h$, apply the first assertion to
$C_{rh}(u,v)$, which divides $C_r(a,b)$ by
Lemma~\ref{primary:cyclotomic-facts}, and use
$\log\max(u,v)=(\log B)/h$. The equality case gives the
factor $h$ when $\rad(h)\mid r$.
For the Fermat assertion, put $t=\tau(a)$, $k=k(a)$,
and write $a=u^{2^t k}$. Then
\[
G_n=\Phi_{2^{n+t+1}k}(u)\mid F_{a,n},\qquad
\varphi(2^{n+t+1}k)=2^{n+t}\varphi(k).
\]
For large $n$, the index is divisible by eight. Applying
the first assertion with $b=1$ to $G_n$ gives the coefficient
$2^t\varphi(k)^2/k$. For each odd prime power $p^e$,
$\varphi(p^e)^2/p^e=p^{e-2}(p-1)^2>1$, proving the
comparison with $2^t$.
\end{proof}

\subsection{Radical bounds}

For integers \(n\geq1\), write \(R_{a,n}=\log\rad(F_{a,n})\) and
\(W_{a,n}=F_{a,n}/\rad(F_{a,n})\).
For every odd prime divisor of \(F_{a,n}\), the order is \(2^{n+1}\),
and lifting the exponent gives
\[
v_p(F_{a,n})=v_p(a^{p-1}-1).
\]
For each such \(p\) and every integer \(k\geq1\), this also gives
\(p^k\mid F_{a,n}\) if and only if
\(a^{p-1}\equiv1\pmod{p^k}\). Since
\(v_2(F_{a,n})\leq1\), for bases of either parity,
\[
\log W_{a,n}
=\sum_{\substack{p\mid F_{a,n}\\p\text{ odd}}}
\{v_p(a^{p-1}-1)-1\}\log p.
\]
The order and valuation identities follow from
Lemma~\ref{primary:cyclotomic-facts}: write \(p-1=2^{n+1}k\),
where \(p\nmid k\), and apply lifting the exponent to
\((a^{2^{n+1}})^k-1\), noting that an odd divisor of
\(a^{2^n}+1\) does not divide \(a^{2^n}-1\).

\begin{proposition}\label{thm:rad-general}\label{thm:fermat-rad}
For every even integer \(a\geq2\) and \(n\geq1\),
\[
P(F_{a,n})\log P(F_{a,n})>2^{n+1}R_{a,n},\qquad
P(F_{a,n})>\frac{2^{n+1}R_{a,n}}{\log(2^{n+1}R_{a,n})}.
\]
For the classical base \(a=2\) and \(n\geq2\), both inequalities
remain valid with \(2^{n+2}\) in place of \(2^{n+1}\).
\end{proposition}

\begin{proof}
Put \(X=P(F_{a,n})\). All prime divisors are \(1\pmod{2^{n+1}}\),
so \(R_{a,n}\leq\lfloor(X-1)/2^{n+1}\rfloor\log X
<X\log X/2^{n+1}\). The classical refinement follows from
Lemma~\ref{lem:quadratic}. In either case, with \(Y=2^{n+1}R_{a,n}\)
or \(Y=2^{n+2}R_{a,n}\), we have \(Y>e\) and
\((Y/\log Y)\log(Y/\log Y)<Y\); monotonicity of \(x\log x\)
gives the second inequality.
\end{proof}

\begin{theorem}\label{thm:radical}\label{thm:rad-asymptotic}
Let \(a\geq2\) be fixed. For every \(\varepsilon>0\)
there are \(0<\gamma=\gamma(\varepsilon)<1\) and
\(r_0=r_0(a,\varepsilon)\) such that, for every
\(r\geq r_0\) with \(8\mid r\), writing
\(q=\varphi(r)\) and \(R=\log\rad\Phi_r(a)\),
\[
R\geq q^{1-\gamma}\quad\Longrightarrow\quad
P(\Phi_r(a))\geq\left(\frac58-\varepsilon\right)qR.
\]
In particular, for all sufficiently large \(n\),
\[
R_{a,n}\geq2^{(1-\gamma)n}\quad\Longrightarrow\quad
P(F_{a,n})\geq\left(\frac58-\varepsilon\right)2^nR_{a,n}.
\]
\end{theorem}
\begin{proof}
It suffices to take \(0<\varepsilon<5/8\). Set
\(Y=(5/8-\varepsilon)qR\), and suppose that
\(P(\Phi_r(a))\leq Y\). Choose \(\rho>0\) so that
\((8/5+\rho/2)(5/8-\varepsilon)<1\), and take \(\delta>0\)
from Lemma~\ref{prop:BT} with tolerance \(\rho\).
Choose \(0<\gamma<1\) with
\(1/(2-\gamma)<1/2+\delta\). The hypothesis and
Lemma~\ref{primary:cyclotomic-facts} give
\[
Y\gg_\varepsilon q^{2-\gamma},\qquad Y\ll_a q^2.
\]
For the fixed odd squarefree part \(b\) of \(a\), put
\(Q=\operatorname{lcm}(2r,b)\). Since
\(r\leq Q\leq2br\) and \(r/q\ll\log\log(3r)\),
\[
Y^{1/2-\delta}\leq Q\leq Y^{1/2+\delta}
\]
once \(r\) is sufficiently large. Lemma~\ref{lem:quadratic}
supplies \(\varphi(Q)/(2q)\) admissible reduced classes.
Summing Lemma~\ref{prop:BT} over them, and allowing for
the possible exceptional prime, gives
\[
R\leq\log r+
\sum_{t\text{ admissible}}\vartheta(Y;Q,t)
\leq\log r+(8/5+\rho/2)Y/q<R.
\]
The last inequality holds for large \(r\), since
\(\log r=o(R)\) under the radical hypothesis.
This contradiction proves the result. Taking
\(r=2^{n+1}\) gives \(q=2^n\) and the Fermat assertion.
\end{proof}

\begin{corollary}\label{cor:eta-general}\label{cor:perfect-power-base}
Let \(a\geq2\) be fixed, and let \(0\leq\eta<1\).
As \(n\to\infty\) through indices satisfying
\(W_{a,n}\leq F_{a,n}^{\eta}\),
\[
P(F_{a,n})\geq
\left(\frac{5\cdot2^{\tau(a)}\log a}{8}
\max_{d\mid k(a)}\varphi(d)
\left(\frac{\varphi(d)}d-\eta\right)
-o_{a,\eta}(1)\right)4^n.
\]
In particular, if \(F_{a,n}\) is squarefree, then
\[
P(F_{a,n})\geq
\left(\frac{5\kappa(a)\log a}{8}-o_a(1)\right)4^n.
\]
\end{corollary}

\begin{proof}
Put \(t=\tau(a)\) and \(k=k(a)\). Choose a divisor
\(d\mid k\) attaining the maximum in the first bound.
The choice \(d=1\) contributes \(1-\eta>0\), so the
chosen divisor satisfies \(\varphi(d)/d>\eta\).
Write \(a=u^{2^t d}\), and set
\[
G_n=\Phi_{2^{n+t+1}d}(u),\qquad
q_n=2^{n+t}\varphi(d).
\]
By Lemma~\ref{primary:cyclotomic-facts}, \(G_n\mid F_{a,n}\).
Thus \(G_n/\rad G_n\leq W_{a,n}\), and
\[
\begin{split}
\log\rad G_n
&\geq\log G_n-\eta\log F_{a,n}\\
&=2^n\left(\frac{\varphi(d)}d-\eta\right)\log a+o_a(2^n).
\end{split}
\]
The right-hand side is a positive fixed multiple of \(2^n\)
up to an \(o_a(2^n)\) term. Since \(q_n\) is also a fixed
multiple of \(2^n\), the radical hypothesis in
Theorem~\ref{thm:radical} holds for large \(n\). Its index
is divisible by eight for large \(n\), and its base \(u\)
is fixed. Therefore
\[
P(F_{a,n})\geq P(G_n)\geq
\left(\frac58-o_{a,\eta}(1)\right)
q_n\log\rad G_n,
\]
which gives the first assertion. There are only finitely many
divisors of \(k\), so their error terms are uniform for fixed
\(a\) and \(\eta\).

If \(F_{a,n}\) is squarefree, then \(W_{a,n}=1\), so take
\(\eta=0\). For every odd prime power \(p^e\),
\(\varphi(p^e)^2/p^e=p^{e-2}(p-1)^2>1\), and this quantity
increases with \(e\). Multiplicativity therefore shows that
\(\max_{d\mid k}\varphi(d)^2/d=\varphi(k)^2/k\), giving
the coefficient \(\kappa(a)\).
\end{proof}

% Unified multiplicity profile and its polynomial benchmarks.
% Mass concentration and bounded-multiplicity consequences belong to thm:mass.

\subsection{Multiplicities when the prime factors are small}
\label{sec:multiplicity}

For an integer $N>1$, put $E(N)=\max_{p\mid N}v_p(N)$.
Define
\[
e_r=\begin{cases}
\log2,&r=2^j,\ j\geq2,\\
\log r,&\text{otherwise}.
\end{cases}
\qquad
\mathcal B_r(Z)=\log Z\min\left\{
\frac Zr,\frac{2Z}{\varphi(r)\log(Z/r)}\right\}\quad(Z>r).
\]
The exceptional prime contributes at most $e_r$ to
$\log \Phi_r(a)$. In particular, for $r=2^{n+1}$ with $n\geq1$,
one has $\Phi_r(a)=a^{2^n}+1$ and the bound is $\log2$, for both
even and odd bases.

\begin{proposition}\label{primary:multiplicity-profile}
\label{primary:cyclotomic-multiplicity-profile}
Let $a\geq2$, $r\geq3$, and suppose that $P(\Phi_r(a))\leq Z$,
where $Z>r$. For every real $R\geq0$,
\[
\sum_{\substack{p\mid \Phi_r(a)\\v_p(\Phi_r(a))\leq R}}
v_p(\Phi_r(a))\log p\leq e_r+R\mathcal B_r(Z).
\]
Consequently,
\[
E(\Phi_r(a))\geq R_*:=
\frac{\log \Phi_r(a)-e_r}{\mathcal B_r(Z)}.
\]
If $R_*>0$, then for every $0<\lambda<1$,
\begin{equation}\label{eq:multiplicity-tail}
\prod_{\substack{p\mid \Phi_r(a)\\v_p(\Phi_r(a))>\lambda R_*}}
p^{v_p(\Phi_r(a))}
\geq\exp\bigl((1-\lambda)(\log \Phi_r(a)-e_r)\bigr).
\end{equation}
\end{proposition}

\begin{proof}
Every nonexceptional prime divisor is $1\pmod r$.
Counting integers congruent to $1$ modulo $r$ and using
Lemma~\ref{lem:classical-BT} give
\[
\vartheta(Z;r,1)\leq\mathcal B_r(Z).
\]
Nonexceptional primes of multiplicity at most $R$ contribute
at most $R\vartheta(Z;r,1)$, and the possible exceptional
prime has valuation one and contributes at most $e_r$.
This proves the first assertion. Take $R=E(\Phi_r(a))$ to obtain
the second. Finally take $R=\lambda R_*$ and subtract the
first bound from $\log \Phi_r(a)$ to prove
\eqref{eq:multiplicity-tail}.
\end{proof}

\begin{corollary}\label{primary:polynomial-multiplicity}
\label{primary:cyclotomic-polynomial-multiplicity}
Fix $\beta\geq0$ and
$(1+\beta)/2\leq\alpha<1+\beta$.
For sufficiently large $A$, if $A<a\leq2A$, $r\geq A^\alpha$,
and $P(\Phi_r(a))\leq A^{1+\beta}$, then
\[
E(\Phi_r(a))\gg_{\alpha,\beta}
\frac{A^{2\alpha-1-\beta}\log A}{(\log\log A)^2}.
\]
If $A<a\leq2A$, $n\geq1$, $2^n\geq A^\alpha$, and
$P(F_{a,n})\leq A^{1+\beta}$ for $F_{a,n}=a^{2^n}+1$,
the stronger bound is
\[
E(F_{a,n})\gg_{\alpha,\beta}A^{2\alpha-1-\beta}\log A.
\]
Under the preceding hypotheses for $\Phi_r(a)$, for each
fixed $0<\lambda<1$ there is $c=c(\alpha,\beta,\lambda)>0$ such that
\[
\prod_{\substack{p\mid\Phi_r(a)\\
v_p(\Phi_r(a))>cA^{2\alpha-1-\beta}\log A/(\log\log A)^2}}
p^{v_p(\Phi_r(a))}
\geq\exp\bigl((1-\lambda)(\log\Phi_r(a)-e_r)\bigr).
\]
Under the preceding hypotheses for $F_{a,n}$, the denominator
$(\log\log A)^2$ in this valuation threshold can be omitted,
and $e_r=\log2$.
\end{corollary}

\begin{proof}
Put $Z=A^{1+\beta}$ and $q=\varphi(r)$.
The cyclotomic size and totient bounds give, uniformly for $r\geq A^\alpha$,
\[
\log \Phi_r(a)-e_r\geq\tfrac12q\log A>0
\]
once $A$ is sufficiently large. Thus $\Phi_r(a)$ has a
nonexceptional prime divisor, which exceeds $r$; the hypothesis
therefore forces $r<Z$. In particular,
$\log\log(3r)\ll_\beta\log\log A$.

The two terms in $\mathcal B_r(Z)$ give
\[
\begin{aligned}
R_*&\geq\frac{q}{2(1+\beta)Z}
\max\{r,\tfrac12q\log(Z/r)\}\\
&\geq\frac{q^2}{6(1+\beta)Z}\bigl(1+\log(Z/r)\bigr).
\end{aligned}
\]
Since $x^2(1+\log(Z/x))$ is increasing for $0<x\leq Z$,
the totient bound and $r\geq A^\alpha$ yield
\[
R_*\gg_\beta
\frac{r^2(1+\log(Z/r))}{Z(\log\log A)^2}
\gg_{\alpha,\beta}
\frac{A^{2\alpha-1-\beta}\log A}{(\log\log A)^2}.
\]
For the Fermat values take $r=2^{n+1}$ and $q=r/2$;
the same calculation omits the factor $(\log\log A)^2$.
The proposition proves the lower bounds for $E$.
Choosing $c$ no larger than $\lambda$ times the respective
lower-bound constant for $R_*$, and using
\eqref{eq:multiplicity-tail}, proves the product lower bounds.
\end{proof}

\section{Cyclotomic values on average}\label{sec:prod}\label{sec:rankavg}

For $A>1$ and $1/2\leq B\leq A$, put
\[
\mathcal S_{A,B}=\{(a,b)\in\mathbb Z^2:A<a\leq2A,
\ B<b\leq2B,\ (a,b)=1\},
\]
and, for $r\geq3$ and $Z\geq2$, define
\[
E_{A,B,r}(Z)=\#\{(a,b)\in\mathcal S_{A,B}:
P(C_r(a,b))\leq Z\}.
\]
The pairs in $\mathcal S_{A,B}$ are ordered. For $B=1/2$,
one has $b=1$. Put
\[
\mathcal S_A=\mathcal S_{A,A},\qquad
E_r^{(1)}(A,Z)=E_{A,1/2,r}(Z),\qquad
E_r^{(2)}(A,Z)=E_{A,A,r}(Z).
\]
M\"obius inversion gives
$\#\mathcal S_A=6A^2/\pi^2+O(A\log A)$.

\subsection{Uniform counting estimates}

\begin{theorem}\label{thm:rectavg}\label{thm:homavg}
For all sufficiently large $A$, uniformly for
$1/2\leq B\leq A$, every integer $r\geq3$, and every $Z\geq2$,
\begin{equation}\label{eq:rectcount}
E_{A,B,r}(Z)\ll\frac{AB}{\varphi(r)}+
\sqrt{\frac r{\varphi(r)}}\,\pi(Z;r,1).
\end{equation}
The first term is omitted when
\begin{equation}\label{eq:rectsharp}
\varphi(r)\geq4\sqrt{AB/\log A}.
\end{equation}
In particular,
\begin{align}
E_r^{(1)}(A,Z)&\ll\frac A{\varphi(r)}+
\sqrt{\frac r{\varphi(r)}}\,\pi(Z;r,1),\label{eq:allcount}\\
E_r^{(2)}(A,Z)&\ll\frac {A^2}{\varphi(r)}+
\sqrt{\frac r{\varphi(r)}}\,\pi(Z;r,1).\label{eq:homcount}
\end{align}
In all these estimates,
\begin{equation}\label{eq:cycpi}
\pi(Z;r,1)\ll\min\left\{\frac Zr,
\frac{Z}{\varphi(r)\log(2+Z/r)}\right\}.
\end{equation}
All constants are absolute, and there is no upper bound
on $r$ or $Z$.
\end{theorem}

\begin{theorem}\label{thm:cycavg}\label{thm:rankavg}\label{thm:allavg}
For every fixed $\beta>2\log2$, choose $0<\eta=\eta(\beta)<1$
such that $\lambda=(1+\eta)\log2<\beta/2$, and put
$c_\eta=(1+\eta)\log(1+\eta)-\eta>0$.
For sufficiently large $A$, depending on $\beta$,
\begin{equation}\label{eq:rankcount}
E_r^{(1)}(A,Z)\ll_\beta\frac A{\varphi(r)}+\pi(Z;r,1)
+\frac A{(\log A)^{c_\eta}}
\qquad\bigl(r\geq(\log A)^\beta\bigr).
\end{equation}
There is also an absolute $A_1$ such that for $A\geq A_1$ and
$r\geq\exp(10\log A/\log\log A)$,
\begin{equation}\label{eq:cyccount}
E_r^{(1)}(A,Z)\ll\frac A{\varphi(r)}+\pi(Z;r,1),
\end{equation}
with the first term omitted for $r\geq\sqrt A$.
In this range,
\begin{equation}\label{eq:cyccut}
E_r^{(1)}(A,Z)\ll\frac A{\varphi(r)}+
\frac{Z}{\varphi(r)\log(2+Z/r)},
\end{equation}
again with the first term omitted for $r\geq\sqrt A$.
The constants in \eqref{eq:cyccount} and \eqref{eq:cyccut}
are absolute. These estimates hold for every $Z\geq2$,
without an upper bound on $r$ or $Z$.
\end{theorem}

\begin{lemma}\label{lem:romega}
Put $T=\log(2A)$, and let $\Omega(a)$ count prime factors with
multiplicity. For fixed $0<\eta<1$,
with $c_\eta=(1+\eta)\log(1+\eta)-\eta$,
\[
\#\{A<a\leq2A:\Omega(a)>(1+\eta)\log T\}
\ll_\eta A/T^{c_\eta}.
\]
For sufficiently large $A$, the ordered divisor function
$d_2$ satisfies
\[
\max_{u\leq(2A)^2}d_2(u)\leq\exp(2T/\log T).
\]
\end{lemma}
\begin{proof}
For $1<z<2$, define a multiplicative function by
$g_z(1)=1$ and $g_z(p^\nu)=(z-1)z^{\nu-1}$ for $\nu\geq1$.
Then $z^{\Omega(n)}=\sum_{d\mid n}g_z(d)$ and $g_z\geq0$, so
\[
\sum_{n\leq X}z^{\Omega(n)}
\leq X\prod_{p\leq X}\left(1+\frac{z-1}{p-z}\right)
\ll_z X(\log X)^{z-1}.
\]
The logarithm of the product is
$(z-1)\sum_{p\leq X}1/p+O_z(1)$; Mertens' estimate gives
the last bound. Take $X=2A$, $z=1+\eta$, and use Markov's
inequality.
For $u=\prod p^{\nu_p}$, put $y=T/(\log T)^3$.
The factors with $p\leq y$ satisfy $\nu_p+1\leq CT$, while
$\nu_p+1\leq2^{\nu_p}$ for $p>y$. Thus
\[
\log d_2(u)\leq y\log(CT)+\frac{2T\log2}{\log y}
=(2\log2+o(1))\frac{T}{\log T}<\frac{2T}{\log T}.
\]
\end{proof}

\begin{lemma}\label{bi-valuation}\label{lem:valuation}
Let $A>1$, $1/2\leq B\leq A$, and $r\geq3$.
Put $q=\varphi(r)$ and $T=\log(4AB)$, and let
\[
N_p(j)=\#\{(a,b)\in\mathcal S_{A,B}:p^j\mid C_r(a,b)\}.
\]
For $p\equiv1\pmod r$ and $p^j>4AB$,
\begin{equation}\label{eq:rectlayer}
N_p(j)\leq\min\left\{q,
2+\frac{8rT^2}{j^2(\log p)^2}\right\}.
\end{equation}
In particular, if $r=2^{n+1}$ with an integer $n\geq1$,
$B=A$, and $T_1=\log(2A)$, then
\begin{equation}\label{eq:layer}
N_p(j)\leq2+\frac{64\,2^nT_1^2}{j^2(\log p)^2}
\qquad(p^j>4A^2).
\end{equation}
If $h=\lfloor\sqrt{2^{n+1}}\rfloor$ and $j\log p>2hT_1$,
then $N_p(j)\leq2$.
For $r=2^{n+1}$ with an integer $n\geq1$, any $1/2\leq B\leq A$ and
$(a,b)\in\mathcal S_{A,B}$, the valuation
$J_p=v_p(a^{2^n}+b^{2^n})$ also satisfies
\begin{equation}\label{eq:bi-max}
J_p\log p\leq\log2+2^n\log(2A).
\end{equation}
\end{lemma}
\begin{proof}
In one root class, two distinct reduced quotients are
separated by at least $p^j/(4B^2)$. They lie in an interval
of length $3A/(2B)$, so
\begin{equation}\label{eq:rectroots}
N_p(j)\leq q(1+6AB/p^j).
\end{equation}
If $p^j>4AB$, Lemma~\ref{lem:fractions} with
$U=2A$ and $V=2B$ gives injectivity of the quotients;
hence $N_p(j)\leq q$.
If $N_p(j)\geq3$, apply Lemma~\ref{lem:growth} with
$s=\lceil j\log p/T\rceil-1\geq j\log p/(2T)$.
The products reduce to the group of $r$th roots, giving
\[
N_p(j)\leq1+\frac{2r}{s(s+1)}
\leq1+\frac{8rT^2}{j^2(\log p)^2}.
\]
For $N_p(j)\leq2$ the asserted bound is immediate.
When $B=A$, $T=2T_1$, which proves \eqref{eq:layer}.
For the assertion $N_p(j)\leq2$, take $s=h$. If $N_p(j)\geq3$,
\eqref{eq:growth} gives at least $(h+1)^2>2^{n+1}$ products,
contradicting the group order.
Finally, $a^{2^n}+b^{2^n}\leq2(2A)^{2^n}$ gives
\eqref{eq:bi-max}.
\end{proof}

\begin{lemma}\label{lem:columns}
Let $A\geq2$, $1/2\leq B\leq A$, $r\geq3$,
$q=\varphi(r)$, $T=\log(4AB)$, and $p\equiv1\pmod r$.
For any $\mathcal T\subseteq\mathcal S_{A,B}$, put
\[
N_p(j)=\#\{(a,b)\in\mathcal T:p^j\mid C_r(a,b)\}.
\]
For every integer $t\geq1$, these counts satisfy
\begin{equation}\label{eq:rectcol}
\sum_{j\geq t}N_p(j)\log p
\ll qABp^{-t}\log p+\sqrt{rq}\,T.
\end{equation}
Let $N_p^{(1)}(j)$ and $N_p^{(2)}(j)$ denote the corresponding
counts for a subset of $(A,2A]\cap\mathbb Z$ and a subset
of $\mathcal S_A$, respectively. With $T_1=\log(2A)$, we have
\begin{equation}\label{eq:plaincol}
\sum_{j\geq t}N_p^{(d)}(j)\log p
\ll qA^dp^{-t}\log p+\sqrt{rq}\,T_1\qquad(d=1,2).
\end{equation}
For a subset $\mathcal G\subseteq(A,2A]\cap\mathbb Z$, suppose
that every product of two members has at most $D^2$ ordered
factorizations with factors in $\mathcal G$, where
$1\leq D\leq\sqrt r$. Its counts also satisfy
\begin{equation}\label{eq:prodcol}
\sum_{j\geq t}N_p(j)\log p
\ll DA\sqrt r\,p^{-t/2}\log p+
(q+\sqrt D\,r^{3/4})T_1.
\end{equation}
All constants are absolute.
\end{lemma}
\begin{proof}
Let $J$ be the largest valuation among the selected pairs,
or zero if the set is empty. Then
$J\log p\leq q\log(4A)\leq2qT$.
Using \eqref{eq:rectroots} when $p^j\leq4AB$ and
\eqref{eq:rectlayer} otherwise gives, for $1\leq j\leq J$,
\begin{equation}\label{eq:alllayers}
N_p(j)\leq10qABp^{-j}+2+
\min\left\{q,\frac{8rT^2}{(j\log p)^2}\right\}.
\end{equation}
The geometric sum handles the first term, and the constant
contributes at most $2J\log p=O(qT)$. Monotonicity gives
\[
(\log p)\sum_{j\geq1}
\min\left\{q,\frac{8rT^2}{(j\log p)^2}\right\}
\leq\int_0^\infty\min\{q,8rT^2/u^2\}\,du
=2T\sqrt{8rq}.
\]
Since $q\leq r$, this proves \eqref{eq:rectcol}.
Taking $B=1/2,A$ proves \eqref{eq:plaincol}.

For \eqref{eq:prodcol}, let $J$ be the largest valuation
among the selected integers, or zero if there are none;
then $J\log p\leq q\log(4A)\leq2qT_1$.
Taking square roots in \eqref{eq:prodcount}, with $s=2$
and $D_2=D^2$, gives
\begin{equation}\label{eq:rankroot}
N_p(j)\leq2DA\sqrt r\,p^{-j/2}+2D\sqrt r.
\end{equation}
If $j\log p>2T_1$, take
$h=\lceil j\log p/T_1\rceil-1\geq j\log p/(2T_1)$
in Lemma~\ref{lem:growth}; this gives
$N_p(j)\leq1+8rT_1^2/(j\log p)^2$.
For $j\log p\leq2T_1$, we have
$8rT_1^2/(j\log p)^2\geq2r\geq2D\sqrt r$.
Combining these bounds with \eqref{eq:rankroot}, for every
$1\leq j\leq J$ we obtain
\[
N_p(j)\leq2DA\sqrt r\,p^{-j/2}+1+
\min\left\{2D\sqrt r,
\frac{8rT_1^2}{(j\log p)^2}\right\}.
\]
The constant costs $O(qT_1)$, and the displayed integral
bounds the minimum by $O(\sqrt D\,r^{3/4}T_1)$.
Summing the geometric term proves \eqref{eq:prodcol}.
\end{proof}

For $r=2^{n+1}$ with an integer $n\geq1$,
Lemma~\ref{lem:columns} with $d=2$ gives
\begin{equation}\label{eq:repeat}
\sum_{j\geq2}N_p^{(2)}(j)\log p
\ll 2^n\log(2A)+A^2\,2^n\frac{\log p}{p^2}
\qquad(p\equiv1\pmod{2^{n+1}}).
\end{equation}

\begin{proposition}\label{prop:moment}\label{bi-moment}
For all sufficiently large $A$, put $q=\varphi(r)$ and
$Y=q^2\log A/100$. For
$(a,b)\in\mathcal S_{A,B}$, define
\[
\mathcal U_{a,b,r}(Z)=
\sum_{\substack{Y<p\leq Z\\p\mid C_r(a,b)}}\log p
+\sum_{\substack{p\leq Z\\p\mid C_r(a,b)}}
(v_p(C_r(a,b))-1)\log p.
\]
Uniformly for all $r\geq3$ and $Z\geq2$,
\begin{equation}\label{eq:rectmoment}
\sum_{(a,b)\in\mathcal S_{A,B}}\mathcal U_{a,b,r}(Z)
\ll\{AB+\sqrt{rq}\,\pi(Z;r,1)\}\log A.
\end{equation}
Every pair with $P(C_r(a,b))\leq Z$ satisfies
\begin{equation}\label{eq:momentlower}
\mathcal U_{a,b,r}(Z)>5q\log A/6.
\end{equation}
For $r=2^{n+1}$ with an integer $n\geq1$, we have
$Y=4^n\log A/100$ and
\begin{equation}\label{eq:moment}
\sum_{(a,b)\in\mathcal S_{A,B}}\mathcal U_{a,b,2^{n+1}}(Z)
\ll\{AB+2^n\pi(Z;2^{n+1},1)\}\log A,
\end{equation}
and every pair with $P(a^{2^n}+b^{2^n})\leq Z$ has
$\mathcal U_{a,b,2^{n+1}}(Z)>5\,2^n\log A/6$.
All constants are absolute, with no upper bound on $r,n,$ or $Z$.
\end{proposition}
\begin{proof}
Put $T=\log(4AB)$, $K=4AB$, and
$Z_0=\min(Z,K)$. The elementary bounds
\begin{equation}\label{eq:elementaryphi}
q\geq\tfrac12r^{3/4},\qquad \log r\leq q
\end{equation}
are useful here. Indeed, $p-1\geq p^{3/4}$ for $p\geq5$;
the factors from $p=2,3$ in $\varphi(r)/r^{3/4}$ have
product at least $2^{1/4}/3^{3/4}>1/2$, and higher
prime powers only increase this ratio. Also
$\log t/t^{3/4}\leq4/(3e)<1/2$ for $t\geq1$.
For $\log A\geq400$, these bounds give $Y\geq r^{3/2}$.
Brun--Titchmarsh consequently yields
\[
\vartheta(u;r,1)\leq6u/q\quad(u\geq Y),\qquad
\vartheta(Y;r,1)\leq3q\log A/50.
\]

For $Y<p\leq Z_0$, use \eqref{eq:rectroots} and partial
summation. If this interval is nonempty, its reciprocal
logarithmic sum is at most
$6(1+\log(Z_0/Y))/q=O(T/q)$, since $Z_0\leq4A^2$.
Hence
\[
\sum_{\substack{Y<p\leq Z_0\\p\equiv1\pmod r}}
N_p(1)\log p\ll\{AB+q\pi(Z_0;r,1)\}T.
\]
For $B=1/2,A$, this also gives
\begin{equation}\label{eq:distinctmass}
\sum_{\substack{Y<p\leq Z_0\\p\equiv1\pmod r}}
N_p^{(d)}(1)\log p
\ll\{A^d+q\pi(Z_0;r,1)\}\log(2A),\qquad d=1,2.
\end{equation}
For integers, the same argument allows $K=(2A)^2$ and
$Z_0=\min(Z,K)$, as used below.

For $p\leq Z_0$ and $j\geq2$, apply
\eqref{eq:rectcol} with $t=2$ and use
\[
\sum_{p\equiv1\pmod r}\frac{\log p}{p^2}
\leq\sum_{l\geq1}\frac{\log(lr+1)}{(lr+1)^2}
\ll\frac{\log(2r)}{r^2}.
\]
Their total is $O(AB+\sqrt{rq}\,\pi(Z_0;r,1)T)$.
For $K<p\leq Z$, use \eqref{eq:rectcol} with $t=1$
to bound $\sum_{j\geq1}N_p(j)\log p$. Since $(\log u)/u$ decreases
for $u>K$,
\[
qAB\frac{\log p}{p}\leq qT/4\leq\sqrt{rq}\,T/4.
\]
Each such prime thus contributes $O(\sqrt{rq}\,T)$.

A prime dividing $r$ has valuation at most one by
Lemma~\ref{primary:cyclotomic-facts}; it also satisfies
$p\leq r<Y$ and contributes nothing to $\mathcal U$.
The preceding three estimates give \eqref{eq:rectmoment},
because $T\asymp\log A$ uniformly in $B$.

For $(a,b)\in\mathcal S_{A,B}$, we also have
\begin{equation}\label{eq:rectsize}
(A/2)^q<C_r(a,b)\leq(4A)^q.
\end{equation}
If $b\leq a/2$, every complex factor has modulus at least
$a/2>A/2$. If $b>a/2$, then
$|a-b\zeta|^2\geq ab|1-\zeta|^2$; multiplication gives
$C_r(a,b)\geq(ab)^{q/2}\Phi_r(1)>(A/\sqrt2)^q$.
For a value with $P(C_r(a,b))\leq Z$, we have
\[
\mathcal U_{a,b,r}(Z)=\log C_r(a,b)
-\sum_{\substack{p\mid C_r(a,b)\\p\leq\min(Y,Z)}}\log p.
\]
The subtracted mass is at most
$3q\log A/50+\log r\leq3q\log A/50+q$. Hence
\[
\mathcal U_{a,b,r}(Z)>
q(\log A-\log2)-3q\log A/50-q>5q\log A/6
\]
for sufficiently large absolute $A$. This proves
\eqref{eq:momentlower}; the Fermat assertions follow by
$q=\varphi(2^{n+1})=2^n$.
\end{proof}

\begin{proof}[Proof of Theorem~\ref{thm:rectavg}]
Summing \eqref{eq:momentlower} over the exceptional pairs
and using \eqref{eq:rectmoment} proves \eqref{eq:rectcount}.

Suppose now that \eqref{eq:rectsharp} holds, and put
$q=\varphi(r)$, $K=AB$, and $T=\log(4AB)$.
Then $q^2\log A\geq16K$.
For $r\leq K^{3/5}$, Brun--Titchmarsh gives
\[
\vartheta(K;r,1)
\leq\frac{2K\log K}{q\log(K/r)}
\leq\frac{5K}{q}\leq\frac5{16}q\log A.
\]
For $K^{3/5}<r\leq K$, progression counting and
\eqref{eq:elementaryphi} give
\[
\frac{\vartheta(K;r,1)}{q\log A}
\leq\frac{K\log K}{rq\log A}
\ll K^{-1/20}=o(1).
\]
For $r>K$ this prime sum is empty. Also,
\[
\frac{\log r}{q\log A}\leq\frac1{\log A},
\qquad
\frac{\sqrt K\log K}{rq\log A}
\leq\frac{\log K}{16\sqrt K}=o(1).
\]
Thus \eqref{eq:lowmass}, including the entire exceptional
prime, gives
\[
\sum_{\substack{p\mid C_r(a,b)\\p\mid r}}\log p
+\sum_{\substack{p^j\mid C_r(a,b),\ p\equiv1\pmod r\\j\geq1,\ p^j\leq K}}\log p
<\tfrac12q\log A
\]
for sufficiently large absolute $A$.

For each $p\equiv1\pmod r$,
\[
\log p\sum_{\substack{j\geq1\\p^j>K}}p^{-j}
\ll T/K.
\]
Indeed, for $p\leq K$ use the first power exceeding $K$;
for $p>K$ use that $(\log u)/u$ decreases for $u>e$.
Summing \eqref{eq:alllayers} over $p^j>K$, and using the
maximum valuation and integral estimate from
Lemma~\ref{lem:columns}, therefore gives
\[
\sum_{\substack{j\geq1\\p^j>K}}N_p(j)\log p
\ll qT+\sqrt{rq}\,T\ll\sqrt{rq}\,T.
\]
For each pair with $P(C_r(a,b))\leq Z$, \eqref{eq:rectsize}
leaves more than $q\log A/3$ in these high layers.
Summing over $p\leq Z$ proves the estimate without its
first term. The prime-counting estimates \eqref{eq:cycpi}
follow from Lemma~\ref{lem:progression}.
\end{proof}

\begin{proof}[Proof of Theorem~\ref{thm:cycavg}]
Put $T=\log(2A)$, $q=\varphi(r)$, $K=(2A)^2$,
$Y=q^2\log A/100$, and $Z_0=\min(Z,K)$.
Let $\mathcal G_A$ and $D$ be given by either
\[
\mathcal G_A=\{A<a\leq2A:\Omega(a)\leq(1+\eta)\log T\},
\quad D=T^\lambda,
\]
or
\[
\mathcal G_A=(A,2A]\cap\mathbb Z,\qquad
D=\exp(T/\log T).
\]
Products of two selected bases have at most $D^2$ ordered
factorizations: use $d_2(u)\leq2^{\Omega(u)}$ in the first
case and Lemma~\ref{lem:romega} in the second case.
In the first range $r\geq(\log A)^\beta$,
$\lambda<\beta/2$ implies, uniformly,
\begin{equation}\label{eq:rankabs}
\frac{\sqrt D\,r^{3/4}}q=o_\beta(1),\qquad
\frac{D\sqrt r}q\log(2+T)=o_\beta(1).
\end{equation}
These ratios are bounded by constant multiples of
\[
T^{\lambda/2}r^{-1/4}\log\log(3r),\qquad
T^\lambda r^{-1/2}\log\log(3r)\log(2+T).
\]
Maximizing the decreasing powers at the lower endpoint proves
both assertions. In the second range, the degree hypothesis
implies $r\geq\exp(8T/\log T)$ for large $A$, and both ratios
again tend uniformly to zero. In either case $D\leq\sqrt r$.

Use \eqref{eq:distinctmass} with $d=1$ and
\eqref{eq:prodcol} with $t=2$ for the repeated logarithms.
Lemma~\ref{lem:rpharm} and \eqref{eq:rankabs} bound their total by
\[
\ll DA\sqrt r\,\frac{T\log(2+T)}q+
(q+\sqrt D\,r^{3/4})\pi(Z_0;r,1)T
\ll\bigl(A+q\pi(Z_0;r,1)\bigr)T.
\]
For $p>K$, $p^{-1/2}\log p\ll T/A$; applying
\eqref{eq:prodcol} with $t=1$ gives
\begin{equation}\label{eq:rankhigh}
\sum_{j\geq1}N_p(j)\log p\ll qT.
\end{equation}
The lower bound \eqref{eq:momentlower}, applied with $b=1$, yields
\begin{equation}\label{eq:rankgood}
\#\{a\in\mathcal G_A:P(\Phi_r(a))\leq Z\}
\ll_\beta A/q+\pi(Z;r,1).
\end{equation}
Lemma~\ref{lem:romega} adds the complement to give
\eqref{eq:rankcount}; the full-interval argument gives
\eqref{eq:cyccount}.

For $r\geq\sqrt A$, take the full interval and replace
$Y$ by $2A$. Then $\vartheta(2A;r,1)=o(q\log A)$ uniformly.
For $r\leq A^{3/4}$, Brun--Titchmarsh gives $O(A/q)$,
whose ratio to $q\log A$ is
\[
O((\log\log A)^2/\log A)=o(1).
\]
For $r>A^{3/4}$, progression counting gives ratio
$O(A\log\log(3r)/r^2)=o(1)$.
Each prime $p>2A$ divides at most $q$ interval bases.
For $p\leq K$, levels $j=2,3$ also contain at most $q$
bases, since $p^2>r^2\geq A$ makes reduction injective
on the interval. Apply \eqref{eq:prodcol} with $t=4$
to the remaining levels. Here $D\leq A^{1/24}$ for
large $A$, $A\leq r^2$, and $\log p\leq2T$, so
\[
DA\sqrt r\,p^{-2}\log p\ll DA r^{-3/2}T
\leq r^{7/12}T=o(qT),\qquad
\sqrt D\,r^{3/4}=o(q).
\]
Thus the sum over $j\geq2$ is $O(qT)$ for every $p\leq K$;
\eqref{eq:rankhigh} handles $p>K$. The unary size estimate and
$\vartheta(2A;r,1)+\log r=o(q\log A)$ leave more than
$q\log A/2$ per counted base, giving
$E_r^{(1)}(A,Z)\ll\pi(Z;r,1)$.
Finally, \eqref{eq:cycpi} proves \eqref{eq:cyccut}.
\end{proof}

An explicit polylogarithmic range is $\beta=3/2$ with
$\eta=1/200$, since $(201/200)\log2<3/4$.

\begin{corollary}\label{cor:cycavg}\label{cor:rankavg}
Let $L(A)\to\infty$ with $\log L(A)=o(\log A)$.
For $r\geq(\log A)^\beta$ as in Theorem~\ref{thm:cycavg},
\[
\#\left\{A<a\leq2A:P(\Phi_r(a))
\leq\frac{A\varphi(r)\log A}{L(A)}\right\}
\ll_\beta\frac{A}{\varphi(r)}+\frac{A}{L(A)}
+\frac{A}{(\log A)^{c_\eta}}=o(A).
\]
For $r\geq\exp(10\log A/\log\log A)$ the last term is
omitted; for $r\geq\sqrt A$ the first is also omitted.
\end{corollary}
\begin{proof}
Put $Z=A\varphi(r)\log A/L(A)$.
If $\varphi(r)/r\leq1/\log A$, progression counting gives
$\pi(Z;r,1)\leq Z/r\leq A/L(A)$.
Otherwise $A/L(A)<Z/r\leq A\log A/L(A)$, so
$\log(2+Z/r)\sim\log A$, and \eqref{eq:cycpi}
gives the same bound. Apply Theorem~\ref{thm:cycavg}.
\end{proof}

\begin{corollary}\label{cor:rectall}
Let $r=r(A)\to\infty$ and $1\leq H(A)=o(\log A)$.
Uniformly for $1/2\leq B\leq A$,
\[
P(C_r(a,b))>AB\frac{\varphi(r)^{3/2}}{\sqrt r}\,H(A)
\]
for almost all $(a,b)\in\mathcal S_{A,B}$.
More precisely, with $q=\varphi(r)$, the exceptional count is
\begin{equation}\label{eq:rectall}
\ll AB/q+ABH(A)/\log A.
\end{equation}
The first term is omitted under \eqref{eq:rectsharp}.
\end{corollary}
\begin{proof}
Put $Z=ABq^{3/2}H/\sqrt r$.
The two bounds in \eqref{eq:cycpi} give
\[
\frac{\sqrt{r/q}\,\pi(Z;r,1)}{AB}
\ll H\min\left\{\frac qr,
\frac1{\log(2+AB(q/r)^{3/2}H)}\right\}
\ll H/\log A.
\]
For $q/r\leq1/\log A$ use the first term; otherwise,
$AB\geq A/2$ and $H\geq1$ make the denominator
comparable to $\log A$. Apply Theorem~\ref{thm:rectavg}.
We have $\#\mathcal S_{A,B}\asymp AB$ uniformly in $B$.
Indeed, counting coprime integers for each denominator gives
\[
\#\mathcal S_{A,B}
=A\sum_{B<b\leq2B}\frac{\varphi(b)}b
+O\bigl(B\log(2+B)\bigr).
\]
The sum is $\asymp B$: M\"obius inversion gives
$6B/\pi^2+O(\log(2+B))$ for large $B$, while bounded
$B\geq1/2$ has at least one denominator and the same
comparison after adjusting the constants.
The error is $o(AB)$ uniformly. Thus \eqref{eq:rectall}
is $o(\#\mathcal S_{A,B})$.
\end{proof}

\subsection{Fermat values}

For every integer $n\geq1$, write
$F_{a,b,n}=a^{2^n}+b^{2^n}=C_{2^{n+1}}(a,b)$ and
$F_{a,n}=F_{a,1,n}$.

\begin{corollary}\label{thm:fermat-avg}\label{thm:almostall}
For all sufficiently large $A$, every integer $n\geq1$,
every $1/2\leq B\leq A$, and every $Z\geq2$,
\begin{equation}\label{eq:fermatrect}
\#\{(a,b)\in\mathcal S_{A,B}:P(F_{a,b,n})\leq Z\}
\ll AB/2^n+\pi(Z;2^{n+1},1).
\end{equation}
The first term is omitted when $2^n\geq4\sqrt{AB/\log A}$.
In particular, the two endpoint estimates are
\begin{align}
\#\{A<a\leq2A:P(F_{a,n})\leq Z\}
&\ll A/2^n+\pi(Z;2^{n+1},1),\label{eq:product-prime-cutoff}\\
\#\{(a,b)\in\mathcal S_A:P(F_{a,b,n})\leq Z\}
&\ll A^2/2^n+\pi(Z;2^{n+1},1).\label{eq:bi-cutoff}
\end{align}
The $A/2^n$ term is omitted when
$2^n\geq4\sqrt{A/\log A}$, and the $A^2/2^n$ term when
$2^n\geq4A/\sqrt{\log A}$.
In all these estimates, $\pi(Z;2^{n+1},1)$ can be replaced by
$Z/[2^n\log(2+Z/2^n)]$.
All constants are absolute, and there is no upper bound on
$n$ or $Z$.
\end{corollary}
\begin{proof}
Take $r=2^{n+1}$ in Theorem~\ref{thm:rectavg}. Since
$\varphi(2^{n+1})=2^n$, its rectangular estimates, including
\eqref{eq:rectsharp}, give the assertions.
The endpoint cases follow by taking $B=1/2$ and $B=A$,
and \eqref{eq:cycpi} gives the stated cutoff estimate.
\end{proof}

In particular, if $2^n\geq4\sqrt{A/\log A}$, then
\begin{equation}\label{eq:product-prime-cutoff-sharp}
\#\{A<a\leq2A:P(F_{a,n})\leq Z\}
\ll\pi(Z;2^{n+1},1)\ll\frac{Z}{2^n\log(2+Z/2^n)}.
\end{equation}
Cutoffs below $2$ give an empty exceptional set.
For fixed $c>0$, uniformly for $1\leq G\leq A$,
\[
\#\{(a,b)\in\mathcal S_{A,B}:P(F_{a,b,n})\leq cAB\,2^n/G\}
\ll_c AB/2^n+\frac{AB}{G\log(2+cAB/G)}.
\]
The first term is absent in the sharper range. For $B=A$,
the denominator is $\asymp_c\log A$, giving the
explicit estimate
\[
\#\{(a,b)\in\mathcal S_A:P(F_{a,b,n})\leq cA^2\,2^n/G\}
\ll_c A^2/2^n+A^2/(G\log A).
\]

\begin{corollary}\label{cor:fermat-all}
Let $L(A)\to\infty$ with $\log L(A)=o(\log A)$.
For every integer $n\geq1$,
\begin{equation}\label{eq:fermatrectall}
\#\{(a,b)\in\mathcal S_{A,B}:
P(F_{a,b,n})\leq AB\,2^n\log A/L(A)\}
\ll AB/2^n+AB/L(A).
\end{equation}
For $B=1/2$ and $B=A$, respectively, this gives
\begin{align*}
\#\{A<a\leq2A:P(F_{a,n})\leq A\,2^n\log A/L(A)\}
&\ll A/2^n+A/L(A),\\
\#\{(a,b)\in\mathcal S_A:
P(F_{a,b,n})\leq A^2\,2^n\log A/L(A)\}
&\ll A^2/2^n+A^2/L(A).
\end{align*}
The first term is absent in the respective sharper range of
Corollary~\ref{thm:fermat-avg}.
If $n=n(A)\to\infty$, all these
exceptional sets have relative density zero, without an
imposed rate of growth or upper bound on $n$.
Under the same assumption $n=n(A)\to\infty$, for every
$1\leq H(A)=o(\log A)$,
\[
P(F_{a,n})>A\,2^nH(A),\qquad
P(F_{a,b,n})>A^2\,2^nH(A)
\]
hold for almost all $a\in(A,2A]\cap\mathbb Z$ and almost all
$(a,b)\in\mathcal S_A$, respectively. The exceptional sets may
depend on $n$.
\end{corollary}
\begin{proof}
At the rectangular cutoff,
$\log(2+Z/2^n)\asymp\log A$ uniformly, since
$A/2\leq AB\leq A^2$ and $\log L=o(\log A)$.
Corollary~\ref{thm:fermat-avg} proves
\eqref{eq:fermatrectall}; the endpoint cutoffs follow
by the same calculation. For the last assertion, take
$L(A)=\log A/H(A)$, giving the explicit rates
\[
\begin{split}
\#\{A<a\leq2A:P(F_{a,n})\leq A\,2^nH(A)\}
&\ll A/2^n+AH(A)/\log A,\\
\#\{(a,b)\in\mathcal S_A:P(F_{a,b,n})\leq A^2\,2^nH(A)\}
&\ll A^2/2^n+A^2H(A)/\log A.
\end{split}
\]
Again, each first term is absent in the respective sharper range.
\end{proof}

\begin{corollary}\label{cor:powers}\label{cor:product-powers}
Fix $\alpha>0$. Uniformly for integers $n\geq1$ with $2^n\geq A^\alpha$
and each fixed $K_0>0$,
\[
\#\{A<a\leq2A:P(F_{a,n})\leq K_0A^{1+\alpha}\}
\ll_{\alpha,K_0}A/\log A.
\]
For each fixed $0<\gamma<\alpha$,
\[
\#\{A<a\leq2A:P(F_{a,n})\leq A^{1+\gamma}\}
\ll_{\alpha,\gamma}A^{1+\gamma-\alpha}/\log A.
\]
If $\alpha\geq1+\gamma$, the latter set is empty for large $A$.
\end{corollary}
\begin{proof}
Use \eqref{eq:product-prime-cutoff} and the monotonicity of
$x/\log(2+x)$. At the second cutoff, if $\alpha<1+\gamma$,
$A^{1-\alpha}$ is absorbed by
$A^{1+\gamma-\alpha}/\log A$.
For $\alpha\geq1+\gamma$, each $F_{a,n}>2$ has an odd
prime divisor exceeding $2^{n+1}$, since its $2$-adic valuation
is at most one. This proves emptiness.
\end{proof}

\subsection{Simultaneous bounds for Fermat values}

\begin{lemma}[{Stewart~\cite[Lemma~3.1, rational case]{StewartLucasLehmer}}]
\label{lem:stewart-rational}
Let $p\geq5$ be prime, let $k\geq1$ be an integer, and let
$x_1,\ldots,x_k$ be multiplicatively independent positive rational
numbers with $v_p(x_i)=0$. For $x=c/d$ in lowest terms with
$c,d>0$, put $h(x)=\log\max(c,d)$. For every integer $N\geq2$,
\[
\begin{split}
v_p\bigl((x_1\cdots x_k)^N-1\bigr)
\leq{}& Ck^3\left(7e\frac{p-1}{p-2}\right)^k
\max\left\{p\left(\frac{k}{\log p}\right)^k,e^k\log p\right\}\\
&\times\prod_{i=1}^k h(x_i)\log N,
\end{split}
\]
where $C$ is absolute.
\end{lemma}

We choose the auxiliary primes by the elementary exchange argument of
Bilu, Gun and Hong~\cite[Lemma~4.1]{BGH}.

\begin{lemma}\label{lem:uniform-stewart}
There is an absolute constant $p_*>0$ such that, for all coprime
integers $u>v\geq1$, every integer $N\geq2$, and every prime
$p\geq p_*$ not dividing $uv$,
\[
v_p(u^N-v^N)<
p\exp\!\left(-\frac{\log p}{51.9\log\log p}\right)
\log u\log N.
\]
\end{lemma}
\begin{proof}
Put
\[
k=\left\lfloor\frac{\log p}{51.8\log\log p}\right\rfloor.
\]
Let $\ell_2,\ldots,\ell_k$ be $k-1$ of the first $k$ primes,
chosen to avoid one fixed prime divisor of $uv$.
For sufficiently large absolute $p$, Chebyshev's elementary
estimate for the $k$th prime gives
\[
\ell_i\ll k\log k<p,\qquad
\log\ell_i\leq1.001\log k\quad(2\leq i\leq k).
\]
The rational numbers
\[
x_1=\frac{u}{v\ell_2\cdots\ell_k},\qquad
x_i=\ell_i\quad(2\leq i\leq k)
\]
are multiplicatively independent $p$-adic units. Indeed,
the chosen prime divisor of $uv$ forces the exponent of $x_1$
in any multiplicative relation to vanish; the $\ell_i$ then
force all the other exponents to vanish.

Since $\prod_{i=1}^k x_i=u/v$ and $p\nmid v$, one has
$v_p((\prod_{i=1}^k x_i)^N-1)=v_p(u^N-v^N)$.
For sufficiently large $p$,
\[
\max\{p(k/\log p)^k,e^k\log p\}=p(k/\log p)^k
\]
so Lemma~\ref{lem:stewart-rational} gives
\[
v_p(u^N-v^N)
\leq Ck^3\left(7e\frac{p-1}{p-2}\right)^k
p\left(\frac{k}{\log p}\right)^k
h(x_1)\prod_{i=2}^kh(x_i)\log N.
\]
The height bounds give
\[
\begin{aligned}
h(x_1)&\leq\log u+1.001(k-1)\log k
\ll(\log u)k\log k,\\
\prod_{i=2}^kh(x_i)&\leq(1.001\log k)^{k-1}.
\end{aligned}
\]
Consequently,
\[
v_p(u^N-v^N)
\leq Ck^4\left(7e\frac{p-1}{p-2}
\frac{1.001k\log k}{\log p}\right)^k
p\log u\log N.
\]
Since
\[
\frac{k\log k}{\log p}\longrightarrow\frac1{51.8},
\qquad \log(7e\cdot1.001/51.8)<-1,
\]
its coefficient is less than
$\exp\{-\log p/(51.9\log\log p)\}$ for all sufficiently
large absolute $p$. All thresholds are independent of $u,v,N$.
\end{proof}

\begin{corollary}\label{cor:tail}\label{cor:ranktail}\label{bi-tail}
For sufficiently large $A$, let $n_0\geq1$ be an integer and define
\[
w_j=(j+1)(1+\log(j+1))^2\qquad(j\geq0).
\]
For every $Z\geq2$, put $Z_j=Z\max\{1,2^j/w_j\}$.
We have
\[
\begin{aligned}
&\#\{A<a\leq2A:\exists j\geq0,\ P(F_{a,n_0+j})\leq Z_j\}
\\[-2pt]
&\hspace{3em}\ll A/2^{n_0}+\frac{Z}{2^{n_0}\log(2+Z/2^{n_0})},
\end{aligned}
\]
\[
\begin{aligned}
&\#\{(a,b)\in\mathcal S_{A,B}:\exists j\geq0,
\ P(F_{a,b,n_0+j})\leq Z_j\}
\\[-2pt]
&\hspace{3em}\ll AB/2^{n_0}+\frac{Z}{2^{n_0}\log(2+Z/2^{n_0})}.
\end{aligned}
\]
The rectangular assertion is uniform for $1/2\leq B\leq A$.
All constants are absolute. The first term is omitted for
$2^{n_0}\geq4\sqrt{A/\log A}$ in the one-variable estimate and
$2^{n_0}\geq4\sqrt{AB/\log A}$ in the rectangle.
Taking $B=A$ gives the estimate over $\mathcal S_A$, with first term
$A^2/2^{n_0}$ and sharper threshold $4A/\sqrt{\log A}$.

Let $L(A)\to\infty$ with $\log L(A)=o(\log A)$.
Outside exceptional sets of sizes $O(A/2^{n_0}+A/L(A))$ and
$O(AB/2^{n_0}+AB/L(A))$, respectively, the inequalities
\begin{align*}
P(F_{a,n})&>
\frac{A\log A}{L(A)}
\max\left\{2^{n_0},\frac{2^n}{w_{n-n_0}},\frac{2^n}{\log A\log\log A}\right\},\\
P(F_{a,b,n})&>
\frac{AB\log A}{L(A)}
\max\left\{2^{n_0},\frac{2^n}{w_{n-n_0}},\frac{2^n}{\log A\log\log A}\right\}
\end{align*}
hold simultaneously for every $n\geq n_0$.
In particular, the lower bounds $A\,2^n/(L(A)\log\log A)$
and $AB\,2^n/(L(A)\log\log A)$ hold throughout the
respective tails.
These exceptional sets have relative density zero whenever
$n_0=n_0(A)\to\infty$. The terms $A/2^{n_0}$ and $AB/2^{n_0}$
are omitted in the respective sharper ranges above.
Since $Z_j\geq Z$, the same counting estimates hold with
every $Z_j$ replaced by $Z$.
\end{corollary}
\begin{proof}
Put $f(x)=x/\log(2+x)$. This function is increasing.
Since $(2+x/2)^3\geq(2+x)^2$, one has
$f(x/2)\leq3f(x)/4$, and consequently
\begin{equation}\label{eq:tail-geometric}
\sum_{j\geq0}f(x2^{-j})\leq4f(x).
\end{equation}
We also claim that
\begin{equation}\label{eq:tail-weighted}
\sum_{j\geq0}f(x/w_j)\ll f(x).
\end{equation}
For bounded $x$, this follows from
$f(x/w_j)\leq x/(w_j\log2)$ and
$\sum_jw_j^{-1}<\infty$.
For sufficiently large $x$, split the sum at $j=\sqrt x$.
When $j\leq\sqrt x$, one has
$w_j\ll\sqrt x(\log x)^2$, so
$\log(2+x/w_j)\gg\log x$.
These terms sum to $O(x/\log x)$.
For the remaining terms, integral comparison gives
\[
\sum_{j>\sqrt x}\frac1{w_j}
\ll\int_{\sqrt x}^\infty
\frac{dt}{t(1+\log t)^2}
\ll\frac1{\log x}.
\]
Using $f(x/w_j)\leq x/(w_j\log2)$ again proves
\eqref{eq:tail-weighted}.

For $n=n_0+j$, apply Corollary~\ref{thm:fermat-avg} with
$Z_j=Z\max(1,2^j/w_j)\geq2$.
Put $V=A$ for integers and $V=AB$ for the rectangle.
In either case, the bound is
\[
\ll \frac{V}{2^{n_0+j}}+
f\left(\frac{Z}{2^{n_0}}\max\{2^{-j},w_j^{-1}\}\right).
\]
The first terms sum to $2V/2^{n_0}$.
By monotonicity of $f$, the sum of the remaining terms is
at most
\[
\sum_{j\geq0}f\bigl((Z/2^{n_0})2^{-j}\bigr)
+\sum_{j\geq0}f\bigl((Z/2^{n_0})/w_j\bigr)
\ll f(Z/2^{n_0}),
\]
by \eqref{eq:tail-geometric} and \eqref{eq:tail-weighted}.
The union bound proves the two counting estimates.
If $2^{n_0}$ satisfies the appropriate sharper threshold,
then every $2^{n_0+j}$ also does, so the first term is absent
for every $j\geq0$.

Finally take $Z=V\,2^{n_0}\log A/L(A)$ in each case.
Since $A/2\leq V\leq A^2$ and
$\log L(A)=o(\log A)$,
$f(V\log A/L(A))\asymp V/L(A)$ uniformly.
Moreover,
\[
2^{n_0}\max(1,2^{n-n_0}/w_{n-n_0})
=\max(2^{n_0},2^n/w_{n-n_0}).
\]
This gives the first two terms in the simultaneous inequalities
and their exceptional counts. We have
$\#((A,2A]\cap\mathbb Z)\asymp A$ and
$\#\mathcal S_{A,B}\asymp AB$ uniformly in $B$,
as proved in Corollary~\ref{cor:rectall}.
Thus the stated relative densities tend to zero.
We finally prove the additional term in the maximum.
For a sufficiently large absolute constant $C_0$, put
\[
R(A)=\exp\{C_0\log(2A)\log\log(3A)\}.
\]
For every integer $n$ with $2^n\geq R(A)$, every
$(a,b)\in\mathcal S_{A,B}$ satisfies
\begin{equation}\label{eq:tail-terminal}
P(a^{2^n}+b^{2^n})>A^2\,2^n.
\end{equation}
The one-variable case is obtained by taking $B=1/2$, so $b=1$.
To see this, put $u=\max(a,b)$ and $v=\min(a,b)$.
Coprimality gives $u>v$, and $A<u\leq2A$.
If all prime factors were at most $A^2\,2^n$, each odd prime
factor would be congruent to $1$ modulo $2^{n+1}$, leaving at most
$A^2/2$ such primes, all satisfying $\log p\ll n$ in this range.
Lemma~\ref{lem:uniform-stewart}, applied
with $N=2^{n+1}$, would bound their total logarithmic mass by
\[
CA^4\,2^n(n\log2)^2
\exp\!\left(-\frac{\log(2^{n+1})}{51.9\log\log(2^{n+1})}\right)\log u.
\]
After increasing the absolute $C_0$, this is $o(2^n\log u)$
uniformly for $2^n\geq R(A)$. Indeed the coefficient decreases
for large $n$, and at the lower endpoint $2^n=R(A)$ it is $A^{4-C_0/51.9+o(1)}$.
Since $v_2(u^{2^n}+v^{2^n})\leq1$, the odd-prime mass is at least
$2^n\log u-\log2$, a contradiction.

Take the cutoff $V\,2^n/(L(A)\log\log A)$.
For large $A$ this exceeds $2$. Corollary~\ref{thm:fermat-avg}
and the progression estimate bound the number of exceptions
at this cutoff by
\[
\ll \frac{V}{2^n}+
f\left(\frac{V}{L(A)\log\log A}\right)
\ll\frac{V}{2^n}+\frac{V}{L(A)\log\log A\log A},
\]
uniformly in both cases. There are
$O(\log A\log\log A)$ indices $n=n_0+j$ with $2^n<R(A)$.
Summing over them gives $O(V/2^{n_0}+V/L(A))$ exceptions.
By \eqref{eq:tail-terminal}, there are no exceptions for
$2^n\geq R(A)$. The $V/2^{n_0}$ term is again absent in the stated sharper
ranges. Taking the union with the preceding exceptional set
proves the three-term maximum.
\end{proof}

\section{Distinct primitive prime divisors}\label{sec:matching}

For a given $r$, we seek a choice of primitive prime divisors
$p_a\mid\Phi_r(a)$, for bases $a\in(A,2A]$, with no prime
chosen twice. When such a choice is impossible for the whole
interval, we ask how many bases must be omitted.

\subsection{Matching deficiency}

\begin{theorem}\label{adaptive:matching}\label{matching:main}
There are absolute constants $C>0$ and $A_0>0$ such that,
for $A\geq A_0$ and every $r\geq3$, pairwise distinct primes
$p_a\mid\Phi_r(a)$ with $\operatorname{ord}_{p_a}(a)=r$
can be chosen for all but at most
\[
C\frac{A^2r}{\varphi(r)^2}
\]
integers $a\in(A,2A]$. For every $Z\geq2$, these primes
can be required to exceed $Z$ after excluding at most
\[
C\frac{A^2r}{\varphi(r)^2}+\pi(Z;r,1)
\]
bases. In particular, if $\varphi(r)^2/(Ar)\to\infty$,
distinct primitive prime divisors can be chosen for almost
all bases. For powers of two, it is enough that $r/A\to\infty$.

For every fixed $\alpha>0$, there is $c_\alpha>0$ such that
the number of excluded bases is $O_\alpha(A^{1-c_\alpha})$
uniformly for all $r\geq A^\alpha$. For every fixed $K>4\log2$ and
$0<\sigma<1-4\log2/K$, the number of excluded bases is
$O_{K,\sigma}\bigl(A\exp\{-(\log A)^\sigma\}\bigr)$,
uniformly for
\begin{equation}\label{eq:matching-subpower}
r\geq\exp\!\left(\frac{K\log A}{\log\log A}\right).
\end{equation}
In both assertions, requiring the primes to exceed $Z$
excludes at most $\pi(Z;r,1)$ additional bases.

In these assertions, for every positive function
$h(S)\to\infty$, the selected primes can also satisfy
\begin{equation}\label{eq:matched-components}
p_a^{v_{p_a}(\Phi_r(a))}>
\frac{\varphi(r)^2\log a}{h(\varphi(r)^2\log a)}.
\end{equation}
The counting constants may be chosen independently of $h$;
the threshold for sufficiently large $A$ may also depend on $h$.
\end{theorem}

\begin{lemma*}[{Shkredov~\cite[Theorem~12 and (35)]{Shkredov19}}]
There is an absolute constant $C_*\geq1$ with the following
property. Let $p$ be an odd prime, let
$G\subset\mathbb F_p^\times$ be a multiplicative subgroup
of order $r$, and put
\[
\rho=\max_{1\leq t<p}
\left|\sum_{x\in G}\exp(2\pi i tx/p)\right|.
\]
For every integer $k\geq2$ satisfying $2^{64k}C_*^4\leq r$,
one has
\begin{equation}\label{eq:shk-input}
\begin{aligned}
\left(\frac{\rho}{r}\right)^{2^{k+1}}
&\leq\frac{2^{4k+6}C_*(\log_2r)^4}{r}\\
&\quad+
\frac{16k^2C_*^{k-1}(\log_2r)^{4(k-1)}p}
{r^{(k+3)/2}},
\end{aligned}
\end{equation}
where $\log_2$ denotes the logarithm to base $2$.
\end{lemma*}
\begin{proof}
Write $\widehat G(t)=\sum_{x\in G}\exp(2\pi i tx/p)$.
It is constant on each nonzero $G$-coset, so
\[
r\rho^{2^{k+1}}
\leq\sum_{t\in\mathbb F_p}|\widehat G(t)|^{2^{k+1}}.
\]
The moment estimate bounds the last sum by
\[
\begin{aligned}
&2^{4k+6}C_*(\log_2r)^4r^{2^{k+1}}\\
&\quad+16k^2C_*^{k-1}(\log_2r)^{4(k-1)}p\,
r^{2^{k+1}-(k+1)/2}.
\end{aligned}
\]
Here we used that $x_1+x_2=x_3+x_4$ has at most $r^3$
solutions in $G^4$. Dividing by $r^{2^{k+1}+1}$ proves the assertion.
\end{proof}

The maximal-valuation selection below is based on the argument of
Erd\H{o}s and Selfridge~\cite[proof of Theorem~3]{ES}.

\begin{proof}[Proof of Theorem~\ref{matching:main}]
Write $q=\varphi(r)$, $T=\log(2A)$, and
$N=\lfloor2A\rfloor-\lfloor A\rfloor$.
Choose either the bipartite divisibility graph between these
bases and their prime divisors $p\equiv1\pmod r$, or, for a fixed
positive function $h(S)\to\infty$, the subgraph retaining just
the edges satisfying \eqref{eq:matched-components}.
In the chosen graph put
\[
D_a=\sum_{p:(a,p)\text{ is an edge}}
 v_p(\Phi_r(a))\log p.
\]
For the full graph, Lemma~\ref{primary:cyclotomic-facts} and
$\log r\leq q$ give $D_a\geq q(\log A-2)$.
For the retained graph, Theorem~\ref{thm:mass} and
Lemma~\ref{lem:hom-size} give the common bound
\[
D_a\geq(1-o(1))q\log A\geq\tfrac12q\log A
\]
uniformly for $r\geq3$. Indeed, $q^2\log a\geq4\log A$,
and deleting the possible exceptional prime costs at most
$\log r\leq q$. Here the vanishing error and the threshold
for $A$ may depend on $h$ in the retained graph. Once the
last inequality holds, all counting constants below are
independent of $h$.

Among matchings with the greatest number of edges, choose one
maximizing the product of the full primary components
$p^{v_p(\Phi_r(a))}$ over its matched pairs $(a,p)$.
Let $U$ be the unmatched bases and put $\delta=\#U$.
Suppose that $\delta>0$, and write $\mathcal N(U)$ for its neighboring primes.
Every $p\in\mathcal N(U)$ is matched to a base $b_p\notin U$;
otherwise an edge could be added to the matching.
For each edge $(a,p)$ with $a\in U$, one has
\[
v_p(\Phi_r(b_p))\geq v_p(\Phi_r(a)),
\]
since otherwise replacing $(b_p,p)$ by $(a,p)$ increases the product.
The neighboring primes are matched to distinct bases, so there are
at most $N-\delta\leq N-1\leq A$ of them. Put
\[
N_p(j)=\#\{a\in U:(a,p)\text{ is an edge},\ p^j\mid\Phi_r(a)\}.
\]
These counts retain the full valuations, and
\[
\sum_{a\in U}D_a
=\sum_{p\in\mathcal N(U)}\sum_{j\geq1}N_p(j)\log p.
\]
Whenever $N_p(j)>0$, the matched base $b_p$ is another root
modulo $p^j$.

For $p>A$, reduction on the interval is injective, so $N_p(j)\leq q$.
Put $x=j\log p/T$. If $x>2$, use the additional root $b_p$
and take $\ell=\lceil x\rceil-1$ in Lemma~\ref{lem:growth} to obtain
\[
N_p(j)\leq\frac{2r}{\ell(\ell+1)}
\leq\frac{2r}{x(x-1)}\leq\frac{4r}{x^2}.
\]
For $x\leq2$, the last bound follows from $N_p(j)\leq q\leq r$.
Thus $N_p(j)\leq\min\{N_p(1),4rT^2/(j\log p)^2\}$,
and comparison with the integral gives
\begin{equation}\label{eq:matchcolumn}
\begin{aligned}
\sum_{j\geq1}N_p(j)\log p
&\leq\int_0^\infty\min\{N_p(1),4rT^2/u^2\}\,du\\
&=4T\sqrt{rN_p(1)}\leq4T\sqrt{r\delta}.
\end{aligned}
\end{equation}

For $r\leq2A$, the inequalities $q\leq r$ and $\delta\leq N\leq2A$
give $\delta\leq4A^2r/q^2$.
For $r>2A$, every neighboring prime exceeds $A$.
Summing \eqref{eq:matchcolumn} over at most $A$ primes gives
\[
\delta(1-o(1))q\log A\leq4AT\sqrt{r\delta},
\qquad
\delta\leq(16+o(1))A^2r/q^2.
\]
This proves the first deficiency bound.
All the matched primes have order $r$
by Lemma~\ref{primary:cyclotomic-facts}. Discarding matched
primes at most $Z$ removes at most $\pi(Z;r,1)$ further bases.

We next sharpen the deficiency estimate when the index grows with $A$.
Allow every $r\geq3$, keeping the chosen graph and the sets above.
The mass lower bound is now
\begin{equation}\label{eq:rowmass}
\sum_{a\in U}D_a\geq(1-o(1))\delta q\log A.
\end{equation}
For $p\leq A$, Lemma~\ref{lem:columns} and comparison with the
integers $1+\ell r$ give
\[
\begin{split}
\sum_{\substack{p\in\mathcal N(U)\\p\leq A}}
\sum_{j\geq1}N_p(j)\log p
&\ll qA\sum_{\substack{p\leq A\\p\equiv1\pmod r}}
\frac{\log p}{p}+\sqrt{rq}\,T\pi(A;r,1)\\
&\ll \frac{qAT^2}{r}+A\sqrt{q/r}\,T.
\end{split}
\]
The sum is empty if $r>A$.

Suppose first that $r\leq A^2$, and put $D=\exp(T/\log T)$.
For each neighboring prime $p>A$, let
\[
G_p=\{x\in\mathbb F_p^\times:x^r=1\},\qquad
\rho_p=\max_{h\in\mathbb F_p^\times}\left|\sum_{x\in G_p}e^{2\pi i hx/p}\right|.
\]
Suppose that $\rho_p\ll r\eta$ uniformly for
$A<p\leq(2A)^2$.
Fourier completion for $(A,2A]\cap\mathbb Z$ gives
\[
\begin{aligned}
N_p(j)\leq N_p(1)
&\leq\frac{r(\lfloor2A\rfloor-\lfloor A\rfloor)}p
+O(\rho_p\log p)\\
&\ll r\{A/p+\eta T\}.
\end{aligned}
\]
By \eqref{eq:matchcolumn}, the mass of each such prime is
\[
\ll rT\{\sqrt{A/p}+\sqrt{\eta T}\}.
\]
There are at most $A$ neighboring primes. In increasing
order they satisfy $p_i\geq1+ir$, whence
$\sum p_i^{-1/2}\ll\sqrt{A/r}$. Their total mass is thus
\[
\ll A\sqrt r\,T+ArT^{3/2}\sqrt\eta.
\]

For $p>(2A)^2$, products of two interval bases are smaller than
$p$, so distinct integer products are distinct modulo $p$.
Lemma~\ref{lem:romega} bounds their ordered multiplicities by
$D^2$. All their residues lie in $G_p$, giving
$N_p(1)\leq D\sqrt r$. Equation~\eqref{eq:matchcolumn} gives mass
$O(\sqrt D\,r^{3/4}T)$ per prime. Combining the three prime ranges with \eqref{eq:rowmass},
and using $q\leq r$ and $D\geq1$, gives
\begin{equation}\label{eq:matchdef}
\delta\ll
\frac{AT}{r}+\frac{Ar}{q}
\left(\sqrt{T\eta}+\sqrt D\,r^{-1/4}\right).
\end{equation}
The implied constant depends only on the uniform constant
in the assumed bound for $\rho_p$.

For $r>A^2$, the first deficiency bound gives
$\delta\ll(\log\log A)^2$ uniformly, since
$(\log\log(3r))^2/r$ is eventually decreasing. This is stronger
than both remaining assertions. In their stated ranges, $r\geq D^2$
for sufficiently large $A$, so it remains to consider
$D^2\leq r\leq A^2$.

For all the neighboring primes $A<p\leq(2A)^2$, take the same integer
\[
k=\left\lceil\frac{4T}{\log r}\right\rceil+4.
\]
Then $r^{k/2}\geq pr^2$ and $k\leq2\log T+5$.
The factors involving $C_*$, $k$ and $\log_2r$ on the right-hand
side of \eqref{eq:shk-input} have logarithm
$O((\log T)^2)=o(\log r)$, uniformly.
The size condition therefore holds, and that side is at most
$r^{-1+o(1)}+r^{-7/2+o(1)}\leq r^{-1/2}$. Hence
\begin{equation}\label{eq:matching-shkredov}
\frac{\rho_p}{r}\leq
\exp\left(-\frac{\log r}{2^{k+2}}\right).
\end{equation}

If $r\geq A^\alpha$ for a fixed $\alpha>0$, then
$k=O_\alpha(1)$, and \eqref{eq:matching-shkredov} gives
$\rho_p\ll_\alpha rA^{-\gamma_\alpha}$ for some
$\gamma_\alpha>0$. Apply \eqref{eq:matchdef} with
$\eta=A^{-\gamma_\alpha}$. Since
$r/q\ll\log\log(3r)$ and $D=A^{o(1)}$, it follows that
$\delta\ll_\alpha A^{1-c_\alpha}$ for some $c_\alpha>0$.

Now fix $K>4\log2$ and suppose that \eqref{eq:matching-subpower} holds.
Put $L=\log T$. Then
\[
k\leq\left(\frac4K+o(1)\right)L+5,
\qquad
\frac{\log r}{2^{k+2}}
\geq T^{1-4\log2/K+o(1)}.
\]
Given $0<\sigma<1-4\log2/K$, choose
$\sigma<\tau<1-4\log2/K$.
By \eqref{eq:matching-shkredov}, we may take
$\eta=\exp(-T^\tau)$ in \eqref{eq:matchdef}.
Here $r/q\ll L$ and $D=\exp(T/L)$, so the final term is at most
\[
A\exp\{-(K/4-1/2+o(1))T/L\}.
\]
Since $K>4\log2>2$, all three terms are
$O_{K,\sigma}(A\exp\{-(\log A)^\sigma\})$.
This completes the proof.

\end{proof}

The graph obtained by imposing \eqref{eq:matched-components} on all
interval bases also admits the following quantitative bound without
the subgroup exponential-sum estimate. Put
$q=\varphi(r)$, $T=\log(2A)$, and $D=\exp(T/\log T)$.
For sufficiently large $A$ depending on $h$ and $r\geq4D^2$, the matching
deficiency satisfies
\begin{equation}\label{eq:matching-elementary}
\delta\ll
\frac{DA^{3/2}\log(Ar+1)}{q\log A}
+\frac{\sqrt D\,Ar^{3/4}}q.
\end{equation}
Indeed, use the retained graph in the proof of
Theorem~\ref{matching:main}, with its unmatched set $U$ and counts $N_p(j)$.
By
\eqref{eq:rankroot},
$N_p(j)\leq f_j+Q$, where
$f_j=2DA\sqrt r\,p^{-j/2}$ and $Q=2D\sqrt r\leq r$.
For $x=j\log p/T\geq2$, the additional witness gives
$N_p(j)\leq4r/x^2$. For $x\leq2$, we have $Q\leq4r/x^2$.
Thus all layers satisfy
\[
N_p(j)\leq f_j+\min(Q,4r/x^2).
\]
The geometric sum and the integral comparison give
\[
\sum_{j\geq1}N_p(j)\log p
\ll DA\sqrt r\,\frac{\log p}{\sqrt p}
+\sqrt D\,r^{3/4}T.
\]
There are at most $A$ neighboring primes. Since
$(\log u)/\sqrt u$ is decreasing for $u>e^2$, comparison
with the first $\lfloor A\rfloor$ terms $1+\ell r$ gives
\[
\sum_{p\in\mathcal N(U)}\frac{\log p}{\sqrt p}
\ll\sqrt{A/r}\,\log(Ar+1).
\]
Summing and dividing by $(1-o_h(1))q\log A$ proves
\eqref{eq:matching-elementary}. Since $D=A^{o(1)}$, this
already gives $\delta=o(A)$ for every fixed
$\varepsilon>0$ and $r\geq A^{1/2+\varepsilon}$, with no
upper restriction on $r$; for $r>A^2$ one may use the
first bound in Theorem~\ref{matching:main}.

Let $1\leq H(A)=o(\log A)$ and put $q=\varphi(r)$.
For $Z=AqH$, progression counting gives
$\pi(Z;r,1)\leq Z/r\leq AH/\log A$ if
$q/r\leq1/\log A$. Otherwise $Z/r\geq A/\log A$, and
Brun--Titchmarsh gives
\[
\pi(Z;r,1)\ll\frac{Z}{q\log(2+Z/r)}
\ll\frac{AH(A)}{\log A}.
\]
Theorem~\ref{matching:main} therefore gives pairwise distinct
primitive prime divisors $p_a>AqH(A)$, still satisfying
\eqref{eq:matched-components} for every positive function
$h(S)\to\infty$, after excluding at most
\[
\begin{cases}
O\!\left(A^2r/q^2+AH(A)/\log A\right),&r\geq3,\\[3pt]
O_\alpha\!\left(A^{1-c_\alpha}+AH(A)/\log A\right),
 &r\geq A^\alpha,\\[3pt]
O_{K,\sigma}\!\left(Ae^{-(\log A)^\sigma}+AH(A)/\log A\right),
 &\displaystyle r\geq\exp\!\left(\frac{K\log A}{\log\log A}\right),
\end{cases}
\]
bases, where $\alpha>0$, $K>4\log2$ and
$0<\sigma<1-4\log2/K$ are fixed. The counting constants are
independent of $h$; the threshold for $A$ may also depend on $h$.
The first bound gives almost all bases when $q^2/(Ar)\to\infty$;
the other two bounds give almost all bases throughout their
stated ranges.

For $r=2^{n+1}$ with $n\geq1$, we have
$\varphi(r)=2^n$ and $\Phi_r(a)=a^{2^n}+1$.
The first bound becomes $O(A^2/2^n+AH(A)/\log A)$ with
$p_a>2^nAH(A)$, so it gives almost all bases when $2^n/A\to\infty$.
The same simultaneous choice is available for every fixed
positive-power range $2^n\geq A^\alpha$ and throughout
\eqref{eq:matching-subpower} with $r=2^{n+1}$.

\subsection{Complete matchings}

\begin{proposition}\label{prop:complete}
For each fixed $\varepsilon>0$ and sufficiently large $A$
depending on $\varepsilon$, and for every integer $r\geq3$,
pairwise distinct primitive prime divisors can be chosen for every
integer base $a\in(A,2A]$ when
\[
\varphi(r)\geq\left(\frac2{\sqrt{3s}}+\varepsilon\right)A\sqrt r,
\qquad s=\gcd(r,2).
\]
Consequently, there is an absolute constant $C>0$ such that,
for all sufficiently large $A$, complete matching holds whenever
\[
r\geq C A^2\bigl(\log\log(3A)\bigr)^2.
\]
For an index $r$ that is a power of two,
$r\geq(8/3+\varepsilon)A^2$ suffices. In these
complete-matching ranges, for every $Z\geq2$,
\[
\#\{A<a\leq2A:P(\Phi_r(a))\leq Z\}\leq\pi(Z;r,1).
\]

For every positive function $h(S)\to\infty$, the selected primes
can also satisfy \eqref{eq:matched-components}. The threshold for
sufficiently large $A$ may also depend on $h$.
\end{proposition}
\begin{proof}
Write $q=\varphi(r)$, $T=\log(2A)$, and
$N=\lfloor2A\rfloor-\lfloor A\rfloor$.
Choose either graph and an extremal matching as in the proof of
Theorem~\ref{matching:main}, with its notation $\delta,U,D_a$. In both cases,
$D_a\geq(1-o(1))q\log A$, with the threshold for $A$ depending
also on $h$ when \eqref{eq:matched-components} is imposed.
Put $s=\gcd(r,2)$.
If $\delta>0$, choose $a_0\in U$. Each prime neighboring $a_0$
is matched to a different base $b_p\ne a_0$ with
$v_p(\Phi_r(b_p))\geq v_p(\Phi_r(a_0))$.
Corollary~\ref{lem:pairmass} gives
\[
v_p(\Phi_r(a_0))\log p
\leq2\sqrt{r/(3s)}\,T+\log s.
\]
There are at most $N-1\leq A$ such primes. Thus
\[
(1-o(1))q\log A
\leq A\{2\sqrt{r/(3s)}\,T+\log s\}.
\]
Since $T/\log A\to1$, this contradicts
$q\geq(2/\sqrt{3s}+\varepsilon)A\sqrt r$ for every fixed
$\varepsilon>0$ and sufficiently large $A$. Hence $\delta=0$.

Finally, $q\gg r/\log\log(3r)$ with an absolute constant.
The function $\sqrt r/\log\log(3r)$ is increasing for large
$r$. At $R=C A^2(\log\log(3A))^2$, it is
$(\sqrt C+o(1))A$. Choosing $C$ sufficiently large and then
$A_0$ ensures the complete-matching inequality for every
$r\geq R$. For a power of two, $s=2$ and $q=r/2$, so
$r\geq(8/3+\varepsilon)A^2$ suffices for each fixed
$\varepsilon>0$ and sufficiently large $A$.
If $P(\Phi_r(a))\leq Z$, its matched prime is at most $Z$;
distinctness proves the final counting assertion.

\end{proof}

\begin{proposition}\label{prop:defect}
Fix $0\leq\eta<1$, $0\leq\lambda<5(1-\eta)/16$ and
$\varepsilon>0$, and put $c=1-\eta-16\lambda/5>0$.
For all sufficiently large $A$ and every integer $r\geq3$ satisfying
\[
\frac{c^2\varphi(r)^2}{c\varphi(r)+2r/\gcd(r,4)}
\geq(1+\varepsilon)A,
\]
pairwise distinct primitive prime divisors can be chosen for every
member of any set $\mathcal B\subset(A,2A]\cap\mathbb Z$ satisfying
\[
\Phi_r(a)/\rad\Phi_r(a)\leq\Phi_r(a)^\eta
\qquad(a\in\mathcal B).
\]
The selected primes can satisfy $p_a>\lambda\varphi(r)^2\log A$.
For every positive function $h(S)\to\infty$, they can also satisfy
\eqref{eq:matched-components}, with the threshold for $A$ depending
on $\eta,\lambda,\varepsilon,h$.
For powers of two $r$, the same conclusions hold whenever
\[
r\geq\left(\frac{2(c+1)}{c^2}+\varepsilon\right)A.
\]
For $\eta=\lambda=0$, so that the values on $\mathcal B$
are squarefree, the same conclusions hold for powers of two $r$
whenever
\[
r\geq(16+\varepsilon)\frac{A}{\sqrt{\log A}}.
\]
\end{proposition}

\begin{proof}
Put $q=\varphi(r)$, $g=\gcd(r,4)$ and $T=\log(2A)$.
Fix a positive function $h(S)\to\infty$, and retain the edges
$(a,p)$ satisfying \eqref{eq:matched-components} and
$p>\lambda q^2\log A$ in the primitive-prime divisibility graph.
Theorem~\ref{thm:mass}, Lemma~\ref{lem:hom-size} and the estimate
for $\vartheta$ in Lemma~\ref{lem:modulus-lifting} give
\[
\sum_{p:(a,p)\text{ retained}}\log p
\geq(c-o_h(1))q\log A\qquad(a\in\mathcal B).
\]
Indeed, the radical hypothesis gives
$\log\rad\Phi_r(a)\geq(1-\eta)\log\Phi_r(a)$.
Deleting small full components costs $o_h(q\log A)$;
deleting the possible exceptional prime costs at most $\log r\leq q$.
For fixed $\lambda>0$, the loss from primes at most
$\lambda q^2\log A$ is at most
$(16\lambda/5+o_\lambda(1))q\log A$.
For $\lambda=0$ this loss is absent.

If Hall's condition fails, choose an inclusion-minimal deficient
set $S\subset\mathcal B$. Put $k=|S|$ and $t=|\mathcal N(S)|$.
The positive row sums imply $k\geq2$. Minimality gives $t=k-1$
and at least two neighbors for every prime in $\mathcal N(S)$.
Let $n_p$ be that number.
Choose two distinct neighbors of $p$ and a prime $\ell_0$
at which their valuations differ. Put
\[
X_p=\{(\log a,v_{\ell_0}(a)):a\text{ is a neighbor of }p\}.
\]
If $X_p$ is collinear, its second coordinate is nonconstant
and hence injective on that line. Thus
$n_p\leq1+T/\log2$, and Corollary~\ref{lem:pairmass} gives
\[
n_p\log p\ll\sqrt r\,T^2.
\]

Suppose that $X_p$ is not collinear, so $n_p\geq3$.
All primitive $r$th roots lie in one coset of the subgroup
of squares in the group of $r$th roots modulo $p$.
After opposite residues are identified, this subgroup has $r/g$
classes: it has order $r/\gcd(r,2)$ and contains $-1$ exactly
when $4\mid r$. Products of exactly $j$ neighbors therefore
occupy at most $r/g$ classes.
The map $a\mapsto(\log a,v_{\ell_0}(a))$ sends products
bijectively to sums. If $p>2(2A)^j$, distinct positive products
remain distinct after opposite residues are identified, since
both their difference and their sum have absolute value less than $p$.
Lemma~\ref{lem:planar} then gives
\begin{equation}\label{eq:exactprod}
(n_p-2)\frac{j(j+1)}2+j+1\leq r/g.
\end{equation}
For $p>4A$, take
$j=\lceil(\log p-\log2)/T\rceil-1\geq1$.
For $p\leq4A$, use $\log p\leq\log2+T$ directly.
In both cases,
\[
\log p\leq\log2+
\frac T2\left(1+\sqrt{1+\frac{8r}{g(n_p-2)}}\right).
\]
Since $n_p/(n_p-2)\leq3$, multiplying by $n_p$ and using
$\sqrt{x+y}\leq\sqrt x+\sqrt y$ gives
\[
n_p\log p\leq n_p\log2+
\frac T2\left(n_p+\sqrt{n_p^2+\frac{8rn_p}{g}}\right)
+O(\sqrt r\,T^2).
\]
The collinear case also satisfies this bound.
Since $n_p\leq k$, $k=t+1$ and $t\leq A$, summing and dividing by $k$ gives
\[
\begin{aligned}
\frac1k\sum_{p\in\mathcal N(S)}n_p\log p
&\leq\frac{tT}{2}\left(1+\sqrt{1+\frac{8r}{gk}}\right)
+O(A+\sqrt r\,T^2)\\
&\leq\frac T2\left(A+\sqrt{A^2+\frac{8rA}{g}}\right)
+O(A+\sqrt r\,T^2).
\end{aligned}
\]
The hypothesis implies $q\geq(1+\varepsilon)A/c$ and $r\geq A$.
By \eqref{eq:totient}, uniformly in this range,
\[
\frac{A+\sqrt r\,T^2}{qT}
\ll\frac1T+\frac{T\log\log(3A)}{\sqrt A}=o(1).
\]
The row lower bound therefore gives
\[
(c-o_h(1))q\log A
\leq\frac T2\left(A+\sqrt{A^2+\frac{8rA}{g}}\right)+o(qT).
\]

Write $s=1+\varepsilon$. The numerical hypothesis is equivalent to
\[
2cq\geq sA+\sqrt{s^2A^2+\frac{8srA}{g}}
\geq\sqrt{s}\left(A+\sqrt{A^2+\frac{8rA}{g}}\right).
\]
Since $T/\log A\to1$, this contradicts the preceding inequality.
Hall's condition therefore holds in the retained graph.
For a power of two, $q=r/2$ and $g=4$, and the fraction in the
hypothesis is $c^2r/(2(c+1))$. Applying the result with
$\varepsilon'=c^2\varepsilon/(2(c+1))$ gives the stated threshold.
For the last assertion, take $\eta=\lambda=0$ and let $r$
be a power of two. The preceding conclusion handles $r\geq5A$,
so suppose that
\[
(16+\varepsilon)A/\sqrt{\log A}\leq r\leq5A.
\]
Use the same retained graph. A minimal Hall-deficient set
$S$ has $k$ bases and $t=k-1$ neighboring primes, each with
degree $n_p\geq2$. Here $q=r/2$, and the row mass is
$(1-o_h(1))q\log A$.

The primes $p\leq8A^2$ contribute at most
$\vartheta(8A^2;r,1)$ per row. Brun--Titchmarsh gives,
uniformly in this range,
\[
\vartheta(8A^2;r,1)
\leq\frac{16A^2\log(8A^2)}{q\log(8A^2/r)}
=(32+o(1))A^2/q.
\]
For each $p>8A^2$, construct $X_p$ from its neighbors as above.
If $X_p$ is collinear, its column mass is $O(\sqrt r\,T^2)$.
Otherwise retain the integer
\[
j=\left\lceil\frac{\log p-\log2}{T}\right\rceil-1\geq2.
\]
It satisfies $p>2(2A)^j$ and $\log p\leq\log2+(j+1)T$.
Since $g=4$, \eqref{eq:exactprod} gives
\[
n_p(j+1)\leq\frac r{2j}+2j-\frac2j
\leq\frac r4+2j.
\]
Corollary~\ref{lem:pairmass} gives $j\leq\log p/T\ll\sqrt r$.
Thus, in both cases,
\[
n_p\log p\leq n_p\log2+\frac{rT}{4}+O(\sqrt r\,T^2).
\]
Summing over $p>8A^2$ and dividing by $k$ gives at most
\[
\frac{rT}{4}+O(A+\sqrt r\,T^2)
=\left(\frac12+o(1)\right)q\log A,
\]
uniformly in the stated range, since $A/(rT)+T/\sqrt r=o(1)$.
The row lower bound now implies
\[
1-o_h(1)
\leq\frac{(128+o(1))A^2}{r^2\log A}+\frac12+o_h(1)
\leq\frac{128}{(16+\varepsilon)^2}+\frac12+o_h(1)<1,
\]
a contradiction. Hall's condition holds, proving the last assertion.
\end{proof}

\section{First moments and repeated prime factors}\label{sec:local}

For the mean-value and density estimates below, let $n\geq1$.
We count positive even bases, for which $F_{a,n}=a^{2^n}+1$ is odd.

\subsection{Mean values}

For $Z\geq2$, put
\[
T_{a,n}(Z)=\sum_{\substack{p\leq Z\\p\mid F_{a,n}}}\log p.
\]

\begin{proposition}\label{prop:first-moment-prime-divisor-mass}
If $A,Z\geq2$ and $2^{n+1}\leq Z^{1-\kappa}$ for fixed
$0<\kappa<1$, then
\begin{equation}\label{eq:first-moment-prime-divisor-mass}
\sum_{\substack{A<a\leq2A\\2\mid a}}T_{a,n}(Z)
\ll_\kappa A\mathcal H(2^{n+1},Z)+Z.
\end{equation}
\end{proposition}

For $z\geq2$, define
\[
\nu_{a,n}(z)=\#\{p\leq z:p^2\mid F_{a,n}\},
\qquad
\lambda_n(z)=2^n\sum_{\substack{p\leq z\\p\equiv1\pmod{2^{n+1}}}}p^{-2},
\]
and let $\lambda_n=\lim_{z\to\infty}\lambda_n(z)$.

\begin{theorem}\label{thm:local-repeated-prime-density}
\label{cor:local-squarefreeness}
For $A,z\geq2$,
\begin{equation}\label{eq:exact-average-repeated-primes}
\begin{aligned}
\sum_{\substack{0<a\leq A\\2\mid a}}\nu_{a,n}(z)
&=\frac A2\lambda_n(z)+O\bigl(2^n\pi(z;2^{n+1},1)\bigr),\\
\lambda_n(z)&\leq\lambda_n<\frac{\pi^2}{24}\,2^{-n}.
\end{aligned}
\end{equation}
Both this sum and the number of even bases $0<a\leq A$ with
$\nu_{a,n}(z)\geq1$ are $O(A/2^n+z)$. If, for fixed
$0<\kappa<1$, $2^{n+1}\leq z^{1-\kappa}$, both are
$O_\kappa(A/2^n+z/\log z)$. For fixed $n,z$, the mean of
$\nu_{a,n}(z)$ over the positive even bases tends to $\lambda_n(z)$.
\end{theorem}

To measure multiplicity, put
\begin{equation}\label{eq:truncated-squarefull-mass}
\Delta_{a,n}(A)=
\sum_p\sum_{\substack{j\geq2\\p^j\leq A\\p^j\mid F_{a,n}}}\log p.
\end{equation}
This is the contribution to $\log(F_{a,n}/\rad(F_{a,n}))$
from prime powers not exceeding $A$.

\begin{proposition}\label{prop:avg-truncated-squarefull-mass}
For $A\geq2$,
\begin{equation}\label{eq:average-truncated-squarefull-mass}
\sum_{\substack{0<a\leq A\\2\mid a}}\Delta_{a,n}(A)
\ll A(n+1)/2^n.
\end{equation}
\end{proposition}

\begin{proof}[Proof of the three statements]
For $p\equiv1\pmod{2^{n+1}}$ and $j\geq1$, put
\[
N_U(p^j)=\#\{0<a\leq U:2\mid a,\ p^j\mid F_{a,n}\}.
\]
Lemma~\ref{primary:cyclotomic-facts} gives $2^n$ roots modulo
$p^j$, hence $2^n$ classes modulo $2p^j$ after imposing parity.
Thus
\[
N_U(p^j)=\frac{U2^{n-1}}{p^j}+O(2^n).
\]
For all other primes the count is zero.

Taking $j=1$ and summing with weight $\log p$ over $(A,2A]$ gives
\[
\sum_{\substack{A<a\leq2A\\2\mid a}}T_{a,n}(Z)
\ll A2^n\sum_{\substack{p\leq Z\\p\equiv1\pmod{2^{n+1}}}}
\frac{\log p}{p}+2^n\vartheta(Z;2^{n+1},1).
\]
Lemma~\ref{lem:rpharm} bounds the first term by
$O(A\mathcal H(2^{n+1},Z))$. Brun--Titchmarsh and
$\log(Z/(2^{n+1}))\geq\kappa\log Z$ give
$2^n\vartheta(Z;2^{n+1},1)\ll_\kappa Z$.

Taking $j=2$ gives
\eqref{eq:exact-average-repeated-primes}. Moreover,
\[
\lambda_n\leq 2^n\sum_{\ell\geq1}(1+2^{n+1}\ell)^{-2}
<2^{-n-2}\sum_{\ell\geq1}\ell^{-2}
=\frac{\pi^2}{24}\,2^{-n}.
\]
Use $\pi(z;2^{n+1},1)\leq z/(2^{n+1})$, or Brun--Titchmarsh under the stated
hypothesis, for the two bounds. Every exceptional base contributes
at least one to the sum. Division by $\lfloor A/2\rfloor$ gives
the mean assertion.

For $j\geq2$ and $p^j\leq A$, the error in $N_A(p^j)$ is
absorbed into $O(A2^n/p^j)$. Summing over $j$ gives
\[
\sum_{\substack{0<a\leq A\\2\mid a}}\Delta_{a,n}(A)
\ll 2^n A\sum_{p\equiv1\pmod{2^{n+1}}}\frac{\log p}{p^2}
\ll 2^n A\sum_{\ell\geq1}
\frac{\log(1+2^{n+1}\ell)}{(1+2^{n+1}\ell)^2}
\ll A\frac{n+1}{2^n}.
\]
This proves
\eqref{eq:average-truncated-squarefull-mass}.
\end{proof}

\subsection{Local squarefreeness}

For $z\geq2$, put
\begin{equation}\label{eq:truncated-local-squarefree-density}
\begin{aligned}
\delta_n(z)&=\prod_{\substack{p\leq z\\p\equiv1\pmod{2^{n+1}}}}
\left(1-\frac{2^n}{p^2}\right),\\
\delta_n&=\prod_{p\equiv1\pmod{2^{n+1}}}
\left(1-\frac{2^n}{p^2}\right).
\end{aligned}
\end{equation}

\begin{proposition}\label{prop:local-squarefree-density}
For fixed $n,z$,
\begin{equation}\label{eq:finite-local-squarefree-density}
\lim_{A\to\infty}
\frac{\#\{0<a\leq A:2\mid a,\ p^2\nmid F_{a,n}
\text{ for every }p\leq z\}}{\lfloor A/2\rfloor}
=\delta_n(z).
\end{equation}
The infinite product converges to a positive number and satisfies
\begin{equation}\label{eq:local-squarefree-density-asymptotic}
1-\delta_n=\lambda_n+O(4^{-n}),\qquad
0\leq1-\delta_n\leq\lambda_n<\frac{\pi^2}{24}\,2^{-n}.
\end{equation}
\end{proposition}

\begin{proof}
The Chinese remainder theorem multiplies the complementary
proportions $1-2^n/p^2$ for the primes $p\leq z$ with
$p\equiv1\pmod{2^{n+1}}$; parity is coprime to these moduli and leaves the proportion
unchanged. This proves \eqref{eq:finite-local-squarefree-density}.

We have $\lambda_n<(\pi^2/24)2^{-n}<1$.
The finite products decrease and are bounded below by
$1-\lambda_n$, so their limit is positive. Since $1-x\leq e^{-x}$,
\[
1-\lambda_n\leq\delta_n\leq e^{-\lambda_n},\qquad
\lambda_n-\lambda_n^2/2\leq1-\delta_n\leq\lambda_n.
\]
The estimate $\lambda_n^2\ll4^{-n}$ proves the assertions.
\end{proof}

The density assertion keeps $n,z$ fixed while $A\to\infty$.
Identifying the infinite product with the density of completely
squarefree Fermat values would require an additional estimate
for the prime-square conditions with $p>z$.

\begin{proposition}\label{prop:subgroup-local-squarefreeness}
Fix $0<\alpha<1$. Uniformly as $A\to\infty$ for integers
$r\geq\max(3,A^\alpha)$,
\begin{equation}\label{eq:subgroup-local-squarefreeness}
\#\{A<a\leq2A:p^2\mid\Phi_r(a)
\text{ for some prime }p\leq A^{1+\alpha/8}\}
\ll_\alpha A^{1-\alpha/4}.
\end{equation}
Taking $r=2^{n(A)+1}$ with $2^{n(A)}\geq A^\alpha$ gives
the same bound for even bases whose
$F_{a,n(A)}=a^{2^{n(A)}}+1$ has a prime-square divisor
$p^2$ with $p\leq A^{1+\alpha/8}$.
\end{proposition}

\begin{proof}
If $p^2\mid\Phi_r(a)$, Lemma~\ref{primary:cyclotomic-facts}
excludes the exceptional prime dividing $r$, whose valuation
is one. Thus $p\equiv1\pmod r$, and $a$ is a primitive
$r$th root modulo $p^2$. These residues belong to the
$r$-element subgroup of all $r$th roots modulo $p^2$.
Lemma~\ref{lem:amplification}, with $s=2$, therefore gives,
for every $\varepsilon>0$,
\begin{equation}\label{eq:subgroup-root-count-large-p}
\#\{A<a\leq2A:p^2\mid\Phi_r(a)\}
\ll_\varepsilon A^\varepsilon\sqrt r\,(1+A/p).
\end{equation}
For the primes $p\leq A^{1+\alpha/8}$ with
$p\equiv1\pmod r$, comparison with the integers
$1+\ell r$, $\ell\geq1$, gives
\[
\pi(A^{1+\alpha/8};r,1)\leq\frac{A^{1+\alpha/8}}r,
\qquad
\sum_{\substack{p\leq A^{1+\alpha/8}\\p\equiv1\pmod r}}
\frac1p\ll\frac{\log A}r.
\]
Summing \eqref{eq:subgroup-root-count-large-p}, and taking
$\varepsilon=\alpha/8$, bounds the exceptional count by a
constant multiple of
\[
\frac{A^{\alpha/8}}{A^{\alpha/2}}
\bigl(A^{1+\alpha/8}+A\log A\bigr)
\ll_\alpha A^{1-\alpha/4}.
\]
\end{proof}

\end{document}